\documentclass[11pt]{amsart}
\usepackage{amsfonts}
\title{On  Klein-Maskit Combination Theorem in space I}
\author{Liulan Li,
Ken'ichi Ohshika and
Xiantao Wang
}
\thanks{The research was partly supported by Program for NCET (No. 04-1783),
NSFs of China (No. 10571048) and Hunan (No. 05JJ10001),\\The third author is the corresponding author}
\address{Liulan Li and Xiantao Wang: Department of Mathematics, Hunan Normal University\\
Changsha,  Hunan 410081, People's Republic of China, Ken'ichi Ohshika: Department of Mathematics,
Graduate School of
Science\\
Osaka University, Toyonaka, Osaka 560-0043, Japan}
\email{lanlimail2008@yahoo.com.cn, ohshika@math.wani.osaka-u.ac.jp, and xtwang@hunnu.edu.cn}
\date{}

\newtheorem{thm}{Theorem}[section]
\newtheorem{lem}[thm]{Lemma}
\newtheorem{pro}[thm]{Proposition}
\newtheorem{cor}[thm]{Corollary}
\theoremstyle{definition}
\newtheorem{defi}{Definition}[section]
\newtheorem{rem}{Remark}[section]
\newtheorem{exa}{Example}[section]
\newtheorem{cl}{Claim}
\newcommand{\C}{\mathbb{C}}
\newcommand{\R}{\mathbb{R}}
\newcommand{\hyp}{\mathbb{H}}
\newcommand{\B}{\mathbb{B}}
\newcommand{\beeq}{\begin{equation}}
\newcommand{\eneq}{\end{equation}}

\numberwithin{equation}{section}
\subjclass{ Primary: 30F40; Secondary: 20H10.}

\keywords{ $N$-dimensional discrete M\"obius group;
Geometrically finite M\"obius group; Block; Standard parabolic
region; Parabolic vertex; Conical limit point; Amalgamated
free product; The Klein-Maskit combination theorem.}

\begin{document}
\begin{abstract}In this paper, we generalise the first
Klein-Maskit combination theorem to  discrete  groups of M\"{o}bius transformations in
higher dimensions. As a simple application of the main theorem, some examples will be constructed.
\end{abstract}

\maketitle

\section{Introduction}
In the theory of classical Kleinian groups, there are theorems
called the combination theorems which give methods to generate new
Kleinian groups as amalgamated free products or HNN
extensions of Kleinian groups. The prototype of such theorems is
Klein's combination theorem which can be rephrased as follows in the
modern terms:
\begin{thm}{(Klein \cite{Kle})}
 Let $G_1$ and $G_2\subset PSL_2 \C$
be two finitely generated Kleinian groups with non-empty regions of discontinuity, and let $D_1$ and
$D_2$ be fundamental domains for $G_1$ and $G_2$ of their regions of discontinuity respectively.
Suppose that the interior of $D_2$ contains the frontier and the
exterior of $D_1$ and that the interior of $D_1$ contains the
frontier and the exterior of $D_2$. Then the group $\langle G_1, G_2 \rangle$ generated by $G_1$ and
$G_2$ in $PSL_2\C$ is a Kleinian group isomorphic to $G_1*G_2$ with non-empty region of discontinuity
and $D=D_1\cap D_2$ is a
fundamental domain for the region of discontinuity of $\langle G_1, G_2 \rangle$.
\end{thm}

Fenchel-Nielsen, in \cite{Fen}, gave a generalisation of Klein's theorem to amalgamated free products
and HNN extensions for Fuchsian groups.
In a series of papers, Maskit considered to generalise Klein's theorem for amalgamated free products
and HNN extensions for Kleinian groups
(\cite{Mas}-\cite{Mas5}).
Thurston gave an interpretation of the combination theorem using three-dimensional hyperbolic geometry
and harmonic maps.
 For applications of the combination theorems, we refer the reader to \cite{Abm, Anc, Bea, Fen,
Koe, Mas6, Yam}.

Among these, the first Maskit combination theorem says that under some conditions  two Kleinian groups
$G_1, G_2$ whose intersection $J$ is geometrically finite generate a Kleinian group isomorphic to the
free product of $G_1$ and $G_2$ amalgamated over $J$ and also under the same conditions the resulting
group is geometrically finite if and only if both $G_1$ and $G_2$ are geometrically finite.

The purpose of the present paper is to generalise this first Maskit
combination theorem to discrete groups of M\"{o}bius transformations
of dimension greater than $2$. A first pioneering attempt to
generalise Maskit's combination theorems to higher dimensions was
made by Apanasov \cite{Apa1, Apa}. In particular he showed that
under the same assumptions as Maskit combined with some extra
conditions, one can get a discrete group which is an amalgamated
free product of two discrete groups of $n$-dimensional M\"{o}bius
transformations. In this paper, we shall show that a generalisation
of Maskit's theorem holds in higher dimensions without any such
additional assumptions, imposing only natural ones. Our theorem also
includes the equivalence of geometric finiteness of the given two
groups and that of  the group obtained by the combination. It should
be noted that in this paper, we say that a Kleinian group is
geometrically finite when the $\varepsilon$-neighbourhood of its
convex core has finite volume for some $\varepsilon>0$, and that we
do not assume that it has a finite-sided fundamental polyhedron. For
more details about these Kleinian groups of higher dimensions, we
refer the reader to \cite{Cqy, Sus, Tuk1, Tuk2, Tuk3} and the
references therein.

Our main result (Theorem \ref{main}) and its proof will appear in
\S4.

This is the first of a series in which we shall discuss
generalisations and applications of Klein-Maskit combination theorem
in higher dimensions. A generalisation of the second Klein-Maskit
combination theorem, which corresponds to HNN extensions, to the
case of discrete groups of M\"{o}bius transformations in higher
dimensions and applications of these two combination theorems will
be given in forthcoming papers.


\section{Preliminaries}
\subsection{Basics on M\"{o}bius transformations}
For $n \geq 2$, we denote by $\bar\R^n$ the one-point
compactification of $\R^n$ obtained by adding $\infty$. The group of
orientation-preserving M\"{o}bius transformations of $\bar\R^n$ is
denoted by $M(\bar \R^n)$, with which we endow the compact-open
topology. We regard $\bar\R^n$ as the boundary at infinity of the
hyperbolic $(n+1)$-space $\hyp^{n+1}$ which is  identified
with the  open unit ball bounded by $\bar \R^n$. We denote the union
of $\hyp^{n+1}$ and $\bar\R^n$ endowed with the natural topology by
$\B^{n+1}$. Any M\"{o}bius transformation of $\bar\R^n$ is extended
to a M\"{o}bius transformation of $\B^{n+1}$, which induces an
isometry of $\hyp^{n+1}$. When it is more convenient, we regard
$\hyp^{n+1}$ as the upper half-space of the $(n+1)$-dimensional
Euclidean space and $\R^n$ as $\{(x_1, \dots , x_n, 0)\}$ in
$\R^{n+1}$. A non-trivial element $g\in M(\bar \R^n)$ is called
\begin{enumerate}
\item {\it loxodromic} if it has two fixed points in $\bar \R^n$
and none in $\hyp^{n+1}$;

\item {\it parabolic} if it has only one fixed point in $\bar
\R^n$ and none in $\hyp^{n+1}$;

\item {\it elliptic} if it has a fixed point in $\hyp^{n+1}$.
\end{enumerate}

For a discrete group $G$ of $M(\bar\R^n)$ and a point $z\in \hyp^{n+1}$ or $x \in \bar\R^n$,
the sets $G(z)=\{g(z)|g\in G\}\subset \hyp^{n+1}$ and $G(x)=\{g(x)|g\in G\}\subset \bar\R^n$
are called $G$-orbits of $z$ and $x$ respectively.
If $z'$ lies in the $G$-orbit of $z$, then we say that $z'$ and $z$ are $G$-equivalent.
\subsection{Limit sets, regions of discontinuity and fundamental sets}
The {\it limit set} $\Lambda (G)$ of a discrete group $G\subset M(\bar \R^n)$ is defined as follows:

$$\Lambda(G)=\overline{G(z)}\cap \bar \R^n$$
for some $z\in \hyp^{n+1}$, where the overline denotes the closure
in $\B^{n+1}=\hyp^{n+1} \cup \bar\R^n$ and $G(z)$ the $G$-orbit of
$z$. We call points of $\Lambda(G)$ {\it limit points}. The
complement $\Omega(G)=\bar \R^n \setminus \Lambda (G)$ is called the
{\it region of discontinuity} of $G$.
The following is a well-known fact.

\begin{lem}
\label{lem:limit set}
Let $G$ be a discrete subgroup of $M(\bar\R^n)$.
If $B\subset\bar \R^n$ is a closed and $G$-invariant subset containing at
least two points, then $\Lambda(G)$ is contained in $B$.
\end{lem}

A discrete group $G\subset M(\bar \R^n)$ is said to act discontinuously at a point $x\in\bar \R^n$ if there is a
neighbourhood $U$ of $x$ such that $\{g\in G|g(U)\cap U\neq\emptyset\}$ is a finite set.
The group $G$ acts discontinuously at every point of $\Omega(G)$, and at no point of $\Lambda(G)$.

The complement of the fixed points of elliptic elements in
$\Omega(G)$ is called the {\it free regular set}, and is denoted by
$^\circ\Omega(G)$. When $^\circ\Omega(G)\neq \emptyset$, a {\it
fundamental set} of $G$ is a set which contains one representative
of each orbit $G(y)$ of  $y\in
$$^\circ\Omega(G)$.
It is obvious that $^\circ\Omega(G)\neq\emptyset$ if and only if
$\Omega(G)\neq\emptyset$.

We have the following lemmata for the limit points. These lemmata
in the classical case when $n=2$ can be found in Theorems II.D.2 and
II.D.5 in Maskit \cite{Mas4}.
 Although the
argument is quite parallel, we give their proofs for completeness.

\begin{lem}
\label{source-sink} Let $x$ be a limit point of a discrete subgroup
$G$ in $M(\bar\R^{n})$. Then there are a limit point $y$ of $G$ and
a sequence $\{g_m\}$ of distinct elements of $G$ such that $g_m$
converges to the constant map to  $x$ uniformly on any compact
subset of $\bar \R^{n+1}\setminus\{y\}$.
\end{lem}

\begin{proof}
Since $x$ is a limit point, there are a point $z\in\hyp^{n+1}$ and a
sequence $\{g_m\}$ of distinct elements of $G$ such that
$g_m(z)\rightarrow x$. Regard $\hyp^{n+1}$ as the upper half-space.
Let $(z_1, \cdots, z_n, z_{n+1})$ be the coordinate of
$z$, with $z_{n+1}>0$. Consider the point $z'=(z_1,\cdots, z_n,
-z_{n+1})$  in the lower half-space. The actions of M\"{o}bius
transformations can be extended to the lower half-space conformally.
Then obviously, we have $g_m(z')\rightarrow x$.

By conjugation, we can assume that $G$ acts on $\B^{n+1}$ with
$\mathrm{Int}\B^{n+1}=\hyp^{n+1}$, that $z=\mathbf{0}$,
and that
$\mathrm{Stab}_G(\mathbf{0})=\mathrm{Stab}_G(\infty)=\{id\}$.
Then $z'=\infty$; hence we have $g_m(\infty)\rightarrow x$. By
taking  a subsequence we can make $g^{-1}_m(\infty)$ converge to
some limit point $y$. Since $g_m$ maps the outside of its isometric
sphere onto the interior of that of $g_m^{-1}$, the  radii of the
isometric spheres of $g_m$ and $g^{-1}_m$, which are equal, converge
to $0$ as $m \rightarrow \infty$, and the centre $g_m(\infty)$ of
the isometric sphere of $g^{-1}_m$ converges to $x$. On the other
hand, the centre of the isometric sphere of $g_m$, which is
$g_m^{-1}(\infty)$ converges to $y$. This completes the proof.
\end{proof}

\begin{lem}
\label{converge to limit points} Let $\{g_m\}$ be a sequence of
distinct elements of a discrete group $G\subset M(\bar \R^n)$. Then
there are  a subsequence of $\{g_m\}$ and limit points $x,y$ of $G$,
which may coincide, such that $g_m$ converges to the constant map
$x$ uniformly on any compact subset of $\bar
\R^{n+1}\setminus\{y\}$.
\end{lem}

\begin{proof}
We may assume that $G$ acts on $\B^{n+1}$ with
$\mathrm{Int}\B^{n+1}$ identified with $\hyp^{n+1}$, and that
$\mathrm{Stab}_G(\infty)=\{id\}$. By taking a subsequence if
necessary, we have two limit points $x$ and $y$ such that
$g_m(\infty)\rightarrow x$ and $g^{-1}_m(\infty)\rightarrow y$. The
conclusion now follows from the proof of Lemma \ref{source-sink}.
\end{proof}

We shall use the following term frequently.

\begin{defi}
Let $H$ be a subgroup of a discrete subgroup $G$ of
$M(\bar\R^n)$. An subset $V$ of $\bar \R^n$ is said to be {\it
precisely invariant} under $H$ in $G$ if $h(V)=V$ for all $h \in H$
and $g(V) \cap V =\emptyset$ for all $g \in G-H$.
\end{defi}

For $\Omega(G)$, we have the following proposition: refer to
Proposition II.E.4 in Maskit \cite{Mas4} or Theorem 5.3.12 in
Beardon \cite{Bea}.

\begin{pro}
\label{domain of discontinuity}
Suppose that $\Omega(G)$ is not empty.
Then a point $x \in \bar\R^n$ is contained in $\Omega(G)$ if and only if

\begin{enumerate}
\item the stabiliser $\mathrm{Stab}_G(x)=\{g\in G| g(x)=x\}$ of $x$ in
$G$ is finite, and

\item there is a neighbourhood $U$ of $x$ in $\bar\R^n$ which is precisely
invariant under $\mathrm{Stab}_G(x)$ in $G$.
\end{enumerate}
\end{pro}


\begin{defi}
A fundamental domain for a discrete group $G$ of $M(\bar\R^n)$ with non-empty region of discontinuity
is an open subset $D$
of $\Omega(G)$ satisfying the following.

\begin{enumerate}
\item $D$ is precisely invariant under the trivial subgroup in $G$.

\item For every $z\in\Omega(G)$, there is an element $g\in G$ such that
$g(z)$ is contained in $\bar D$, where $\bar D$ denotes the closure of $D$ in $\bar
\R^n$.

\item $\mathrm{Fr} D$, the frontier of $D$ in $\bar \R^n$, consists of
limit points of $G$, and a finite or countable collection of
codimension-$1$ compact smooth submanifolds with boundary, whose boundary is contained in
$\Omega(G)$ except for a subset with $(n-1)$-dimensional Lebesgue measure $0$.
The intersection of each submanifold with
$\Omega(G)$ is called a side of $D$.

\item For any side $\sigma$ of $D$, there are another side $\sigma'$ of $D$, which may
coincide with $\sigma$, and  a nontrivial element $g\in G$ such that
$g(S)=S'$. Such an element $g$ is called the side-pairing
transformation from $\sigma$ to $\sigma'$.

\item If $\{\sigma_m\}$ is a sequence of distinct sides of $D$, then the
diameter of $\sigma_m$ with respect to the ordinary spherical metric on $\bar \R^n$ goes to $0$.

\item For any compact subset $K$ of $\Omega(G)$, there are only finitely many translates of
$D$ that intersect $K$.
\end{enumerate}
\end{defi}

A fundamental set $F$ for a discrete subgroup $G$ of $M(\bar \R^n)$ whose
interior is a fundamental domain is called {\it a constrained
fundamental set}.

\subsection{Normal forms}
Let $G_1$ and $G_2$ be two subgroups of $M(\bar \R^n)$, and $J$ a subgroup of $G_1 \cap G_2$.

A {\it normal form} is a word consisting of alternate products of
elements of $G_1-J$ and those of $G_2-J$. Two normal forms
$g_n\cdots g_kg_{k-1}\cdots g_1$ and $g_n\cdots
(g_kj)(j^{-1}g_{k-1})\cdots g_1$ are said to be {\it equivalent} for
any $j\in J$. The word length of the normal form is simply
called the {\it length}. The length is invariant under the
equivalence relation.

A normal form is called a $1$-form if the last letter is
contained in $G_1-J$, and a $2$-form otherwise. More
specifically a normal form is called an $(m, k)$-form if the last
letter is contained in $G_m-J$ and the first letter is contained in
$G_k-J$.

The multiplication of two normal forms is defined to be the
concatenation of two words which is contracted to the minimum length by the equivalence defined above.
The product of two normal forms is equivalent to
either a normal form or to an element of $J$.

It is obvious that any element of the free product of $G_1$ and $G_2$ amalgamated over $J$,
which is denoted by $G_1 *_J G_2$, either is an element of $J$ or can be expressed in a normal
form, and that there is a one-to-one correspondence between $G_1 *_J G_2$ and the union of $J$
and the set of the equivalence classes of normal forms.
Also it is easy to see that this correspondence is an isomorphism with respect to the
multiplication defined above.



Let $\langle G_1, G_2 \rangle$ denote the subgroup of $M(\bar\R^n)$ generated by $G_1$ and $G_2$.
There is a natural homomorphism
$\Phi: G_1*_JG_2\rightarrow
\langle G_1, G_2\rangle$
 which is defined by $\Phi(g_n\cdots
g_1)=g_n\circ\cdots \circ g_1$ for a normal form $g_n \dots g_1$
representing an element of $G_1 *_J G_2$, and $\Phi(j)=j$ for $j \in
J$. It is easy to see that this is well defined and
independent of a choice of a representative of the
equivalence class. The map is obviously an epimorphism.

If $\Phi$ is an isomorphism, then we write $\langle G_1,
G_2\rangle=G_1*_JG_2$ identifying elements of $G_1 *_J G_2$ and their images by $\Phi$.

Since $J$ is embedded in $\langle G_1, G_2\rangle$, each nontrivial
element in the kernel of $\Phi$ can be written in a normal
form.

\begin{lem}
\label{kernel}
$\langle G_1, G_2\rangle=G_1*_JG_2$ if and only if $\Phi$ maps no
non-trivial normal forms to the identity.
\end{lem}



\subsection{Interactive pairs}
Following Maskit, we shall define interactive pairs as follows.

Let $G_1$ and $G_2$ be two discrete subgroups of $M(\bar
\R^n)$ and $J$ a subgroup of $G_1 \cap G_2$ as in
the previous subsection. Let $X_1, X_2$ be disjoint non-empty
subsets of $\bar\R^n$. The pair $(X_1, X_2)$ is said to be an {\it
interactive pair} (for $G_1, G_2, J$) when
\begin{enumerate}
\item each  of $X_1, X_2$ is
invariant under $J$,
\item every element of $G_1-J$ sends $X_1$ into $X_2$,
\item
and every element of $G_2-J$ sends $X_2$ into $X_1$.
\end{enumerate}
 An interactive
pair is said to be {\it proper} if there is a point in $X_1$ which is not contained in a $G_2$-orbit
of  any point of $X_2$, or there is a point in $X_2$
which is not contained in a $G_1$-orbit of any point of $X_1$.

\begin{lem}{(Lemma $VII.A.9$ in \cite{Mas4})}
\label{interactive}
Suppose that $(X_1, X_2)$ is an interactive pair for $G_1, G_2, J$.
Let $g=g_n\cdots g_1$
be an $(m, k)$-form.
Then we have $\Phi(g)(X_k)\subset X_{3-m}$.
Furthermore if $(X_1,X_2)$ is proper and $g$ has length greater than $1$, then
the inclusion is proper.
\end{lem}

The existence of a proper interactive pair forces $\Phi$ to be isomorphic.
(Theorem $VII.A.10$ in Maskit \cite{Mas4} in the case when $n=2$.)

\begin{thm}
\label{free product}
Let $G_1,G_2, J$ be as above and suppose that there is a proper interactive pair for $G_1, G_2, J$.
 Then
$\langle G_1, G_2\rangle=G_1*_JG_2$.
\end{thm}

This  easily follows from Lemmata \ref{kernel} and
\ref{interactive}.

The following is a straightforward generalisation of Theorem VII.A.12 in
Maskit \cite{Mas4}.

\begin{lem}
\label{precisely invariant} Suppose  that $(X_1,X_2)$ is
an interactive pair for $G_1, G_2, J$. Suppose moreover that there
is a fundamental set $D_m$ for $G_m$ for $m=1,2$ such that
$G_m(D_m\cap X_{3-m})\subset X_{3-m}$. Then $D=(D_1\cap
X_{2})\cup(D_2\cap X_{1})$ is precisely invariant under $\{id\}$ in
$G=\langle G_1,G_2\rangle$. Furthermore, if $D$ is non-empty, then
$\Phi$ is isomorphic.
\end{lem}

\begin{proof} What we shall show is that for any $x \in D$ and any non-trivial element $g \in G_1*_J G_2$,
we have $\Phi(g)(x) \not\in D$. Since this holds trivially for the
case when $D$ is empty, we assume that $D$ is non-empty. We assume
that $x$ is contained in $D_1 \cap X_2$. The case when $x$ lies in
$D_2 \cap X_1$ can be dealt with in the same way.

If $g$ is a non-trivial element in $J$, then $g(x)$ lies in $X_2$ since
$X_2$ is $J$-invariant.
On the other hand, since $D_1$ is a fundamental set, we have $g(x) \not\in D_1$.
These imply that $g(x) \not\in D$.

Now we shall consider the case when $g$ is represented in a normal form.

\begin{cl}
If $g=g_ng_{n-1}\cdots g_1$ is an $m$-form ($m=1$ or $2$), then
$\Phi(g)(x)\in X_{3-m}\setminus D_m$.
\end{cl}

\begin{proof}We shall prove this claim by induction.

We first consider the case when $n=1$. Suppose first that  $g$ is an
element in $G_1-J$. Then $\Phi(g)(x)\in X_2$ by assumption, whereas
$\Phi(g)(x) \not\in D_1$ since $D_1$ is a fundamental set of
$G_1$. Therefore $\Phi(g)(x)$ is not contained in $D$ in
this case. Suppose next that $g$ is in $G_2-J$. Then $\Phi(g)(x)$
lies in $X_1$ since the assumption that $(X_1,X_2)$ is an
interactive pair implies $\Phi(g)(X_2) \subset X_1$. We shall show
that $\Phi(g)(x)$ does not lie in $D_2$. Suppose, seeking a
contradiction, that $\Phi(g)(x)$ lies in $D_2$. Then since
$\Phi(g^{-1})$ is contained in $G_2-J$ and $\Phi(g)(x) \in X_1 \cap
D_2$, by assumption, we have $x=\Phi(g^{-1})\Phi(g)(x)$ lies in
$X_1$. This contradicts the assumption that $x$ lies in
$X_2$.

Now, we assume that our claim holds in the case when $g$ has length
$n-1$, and suppose that $g$ has length $n$. We consider the case
when $g$ is a $(3-m)$-form. The case when $g$ is an $m$-form can
also be dealt with in the same way. Since $\Phi (g_{n-1}\cdots
g_1)(x)\in X_{3-m}\setminus D_m$ by the assumption of induction, we
have $\Phi(g)(x)\in g_n(X_{3-m}\setminus D_m)\subset X_m$.

Suppose that $\Phi(g)(x)$ lies in $D_{3-m}$.
Then we have $\Phi(g)(x)\in X_m\cap D_{3-m}$.
This implies that
$\Phi(g_{n-1}\cdots g_1)(x)\in g_n^{-1}(X_m\cap D_{3-m})\subset X_m.$
This is a contradiction.
Thus we have shown that $\Phi(g)(x)$ is contained in $X_m \setminus D_{3-m}$.
\end{proof}
By what we have proved above, if $D\neq \emptyset$, then for any
$g\in G_1*_JG_2-\{id\}$, we have $\Phi(g)(D)\cap
D=\emptyset$. This in particular shows that $\Phi(g)\neq id$. Then
Lemma \ref{kernel} shows that $G=G_1*_JG_2$.
\end{proof}

\begin{rem}
Maskit called a fundamental set $D_m$ for $G_m$  maximal with
respective to $X_m$ (which is precisely invariant under
$J$ in $G_m$) if $D_m\cap X_m$ is a fundamental set for
the action of $J$ on $X_m$, and in Theorem VII.A.12 in \cite{Mas4},
the fundamental sets $D_1, D_2$ were assumed to be maximal. Also the
proof of the theorem above shows that the assumption of
maximality is in fact redundant.
\end{rem}

In Maskit \cite{Mas4}, the following sufficient condition for two
open balls to be an interactive pair is given.

\begin{pro}[Proposition VII.A.6 in \cite{Mas4}]
\label{Maskit}
Let $G_m\subset M(\bar\R^n)$ $(m=1,2)$ be two discrete groups with
a common subgroup $J$ and $S\subset \bar\R^n$ be an $(n-1)$-sphere
bounding two open balls $X_1$ and $X_2$. If each $X_m$ is precisely
invariant under $J$ in $G_m$, then $(X_1,X_2)$ is an interactive
pair.
\end{pro}


\subsection{Convex cores and geometric finiteness}

%
%
%
%
%
\begin{defi}
\label{convex core} Let $G$ be a discrete subgroup of $M(\bar\R^n)$
and $\Lambda(G)$ its limit set. We denote by
$\mathrm{Hull}(\Lambda(G))$, the minimal convex set of $\hyp^{n+1}$
containing all geodesics whose endpoints lie on $\Lambda(G)$. This
set is evidently $G$-invariant, and its quotient
$\mathrm{Hull}(G)/G$ is called the {\em convex core} of $G$, and is
denoted by  $\mathrm{Core}(G)$. The group $G$ is said to be {\it
geometrically finite} if there exists
$\varepsilon>0$ such that the
$\varepsilon$-neighbourhood of $\mathrm{Core}(G)$ in
$\hyp^{n+1}/G$ has finite volume.
\end{defi}

As we shall see below, Bowditch proved in \cite{Bow} that this condition is equivalent to
other reasonable definitions of geometric finiteness, except for the one that $\hyp^{n+1}/G$
has a finite-sided fundamental polyhedron, whose equivalence to the above condition has
not been known until now.
%
\subsection{Euclidean isometries}
The classification of discrete groups of Euclidean isometries is known as  Bieberbach's
theorem (see \cite{Wol} or \cite{Rat}, for example).

\begin{thm}{(Bieberbach)}
\label{Bieberbach}
Let $G$ be a discrete group of
Euclidean isometries of $\R^n$.
Then the following hold.
\begin{enumerate}
\item If $\R^n/ G$ is compact, then there is a normal subgroup
$G^{*}\subset G$ of finite index consisting only of
Euclidean translations, which is isomorphic to a free abelian group
of rank $n$.

\item If $\R^n/ G$ is not compact, then there exists a normal
subgroup $G^{*}\subset G$ of finite index in $G$ which is a free
abelian group of rank k with $0\leq k \leq n-1$.
\end{enumerate}

By taking conjugates of $G$ and $G^*$ with respect to an isometry of $\R^n$, the groups
can be made to
have the following properties.

Decompose $\R^n$ into $\R^k \times \R^{n-k}$, where $\R^k$ is identified with  $\R^k\times \{0\}\subset
\R^n$ and $\R^{n-k}$ with $\{0\}\times \R^{n-k}\subset \R^n$.
Let $g(x)= U(x)+ a $ be an
arbitrary element of $G$, where $U$ is a rotation and $a$ is an element of $\R^n$.
Then the rotation $U$ leaves $\R^k$ and
$\R^{n-k}$ invariant and the vector $a$ lies in the subspace $\R^k$.
Furthermore,
 if $g$ lies in $G^{*}$, then $U$ acts on $\R^k$ trivially.
\end{thm}
In the following we always identify the factors of the decomposition
$\R^n =\R^k \times \R^{n-k}$ with $\R^k\times \{0\}$ and $\{0\}\times \R^{n-k}$.

\begin{defi}
\label{maximal abelian}
For a discrete subgroup $G$ of Euclidean isometries, we define $G^*$ to be a free
abelian normal subgroup of $G$ which is maximal among those having the property in Theorem \ref{Bieberbach}.
\end{defi}

\subsection{Extended horoballs, peak domains and standard parabolic regions}
A point $x$ of $\Lambda(G)$ of a discrete group $G$ of
M\"{o}bius transformations is called {\em a parabolic fixed point}
if $\mathrm{Stab}_G(x)$ contains parabolic elements. An easy
argument shows that $\mathrm{Stab}_G(x)$ cannot contain a loxodromic
element then. For a parabolic fixed point $z$, a horoball
in $\B^{n+1}$ touching $\bar \R^n$ at $z$ is invariant
under $\mathrm{Stab}_G(z)$.
In the case when
$\mathrm{Stab}_G(z)$ has rank less than $n$, it is useful
to consider a domain larger than a horoball as follows.

\begin{defi}
Let $G$ be a discrete subgroup of $M(\bar \R^n)$. Let $z$ be a point
of $\bar \R^n$ which is not a loxodromic fixed point. Let
$\mathrm{Stab}^*_G(z)$ be the maximal free abelian subgroup as in
Definition \ref{maximal abelian} of the stabilizer
$\mathrm{Stab}_G(z)$ of $z$ in $G$. Suppose that the rank of
$\mathrm{Stab}^*_G(z)$ is $k$ with $k \leq n-1$. Then there is a
closed subset $B_z\subset \B^{n+1}$ invariant under
$\mathrm{Stab}_G(z)$ which is in the form

$$B_z
=h^{-1}\{x\in \B^{n+1}| \sum_{i=k+1}^{n+1} {x_i}^2\geq
t\},$$ where $t$ $(>0)$ is a constant and $h\in M(\bar \R^n)$ is  a
M\"obius transformation such that $h(z)=\infty$. We call $B_z$  an
{\it extended horoball} of $G$ around $z$.
\end{defi}

\begin{defi}
Let $T_1,\cdots, T_m$ be subsets of $\bar R^n$ and $J_1,\cdots, J_m$
subgroups of the group $G\subset{M(\bar \R^n)}$. We say that
$(T_1,\cdots, T_m)$ is precisely invariant under $(J_1,\cdots, J_m)$
in $G$, if each $T_k$ is precisely invariant under $J_k$ in $G$, and
if for $i\neq j$ and  all $g\in G$, we have  $g(T_i)\cap
T_j=\emptyset.$
\end{defi}

\begin{defi} A {\em peak domain} of a discrete
group $G$ of $M(\bar\R^n)$ with non-empty region of discontinuity at the parabolic fixed
point $z\in \bar \R^n$ is an open
subset $U_z \subset \bar \R^n$ such that
\begin{enumerate}
\item $U_z$ is precisely invariant under $\mathrm{Stab}_G(z)$ in $G$.

\item there exist a $t> 0$, and a transformation $h\in M(\bar \R^n)$ with
$h(z)=\infty$ such that
$$\{x\in \R^n| \sum ^{n}_{i=k+1} {x_i}^{2} >  t \}= h(U_z),$$
where $k=\mathrm{rank}\, \mathrm{Stab}^*_G(z), 1\leq k\leq n-1.$
\end{enumerate}
\end{defi}

\begin{defi}
If $G$ has an extended horoball $B$ around $z$, then the interior of
its intersection with $\bar \R^n$ is a peak domain. Following
Bowditch \cite{Bow}, we use the term {\em standard parabolic region}
at $z$ to mean an extended horoball when the rank of $Stab_G(z)$ is
less than $n$, and a horoball when the rank of $Stab_G(z)$ is $n$.

\end{defi}

\begin{defi}
\label{parabolic vertex}
A point $z\in \bar \R^n$ fixed by a parabolic element of a discrete group
$G \subset M(\bar\R^n)$ is
said to be a {\em parabolic vertex} of $G$ if one of the following conditions is
satisfied.
\begin{enumerate}
\item The subgroup $\mathrm{Stab}^*_G(z)$ has rank $n$.

\item There exists a peak domain $U_z$ at the point $z$.
\end{enumerate}
\end{defi}

\begin{rem}
It is easy to see that the two conditions in Definition
\ref{parabolic vertex} are mutually exclusive: a peak domain exists
only if $\mathrm{rank}\,\mathrm{Stab}^*_G(z)<n$. Also we can easily
see that, in the case when $n=2$, the definition coincides with that
of  cusped parabolic fixed points as in Beardon-Maskit
\cite{Bea2}.
\end{rem}

\begin{defi}
A parabolic fixed point $z$ for the group $G$ is called bounded if
$(\Lambda(G)\setminus\{z\})/\mathrm{Stab}_G(z)$ is compact (see Bowditch
\cite{Bow, Bow1}).
\end{defi}
There is a relationship between a bounded parabolic fixed point and
a parabolic vertex, which was proved by Bowditch \cite{Bow}.

\begin{lem}
\label{bounded} $z$ is a bounded parabolic fixed point for a
discrete group $G$ if and only if $z$ is a parabolic
vertex.
\end{lem}

\begin{defi}
Let $G$ be a discrete subgroup of $M(\bar\R^n)$. A point $x \in
\bar\R^n$ is said to be a conical limit point (or a point of approximation in some literature) if  there are $z\in
\hyp^{n+1}$ and a geodesic ray $l$ in $\hyp^{n+1}$ tending
to $x$ in $\B^{n+1}$ whose $r$-neighbourhood with some $r\in \R$
contains infinitely many translates of $z$.
\end{defi}

Conical limit points can be characterised as follows.
See Theorem 12.2.5 in Ratcliffe \cite{Rat}.

\begin{pro}
\label{conical limit} Let $G$ be a discrete group of $M(\bar \R^n)$
regarded as acting on $\B^{n+1}$ by hyperbolic isometries. Then a
point $z \in \partial \B^{n+1}$ is a conical limit point of $G$ if
and only if there exist $\delta >0$, distinct elements $g_m$ of $G$,
and $x \in \partial \B^{n+1}\setminus\{z\}$ such that
$g_m^{-1}(\mathbf{0})$ converges to $z$ while $|g_m(x)-g_m(z)| >
\delta$ for all $m$. Furthermore, if this condition holds, then for every
$x \in \partial \B^{n+1}\setminus\{z\}$, there is
$\delta>0$ such that $|g_m(x)-g_m(z)| > \delta$ for all $m$.
\end{pro}

The following result due to Bowditch \cite{Bow} or \cite{Bow1} will
be essentially used in the proof of our main theorem.

\begin{pro}
\label{limit set}
Let $G\subset M(\bar \R^n)$ $(n\geq 2)$ be a discrete group.
Then $G$ is geometrically finite if and only if every point of $\Lambda(G)$ is
either a parabolic vertex or a conical limit point.
\end{pro}

\subsection{Dirichlet domains and standard parabolic regions}
Dirichlet domains are fundamental polyhedra of hyperbolic
manifolds, which will turn out to be very useful for us.

\begin{defi}
\label{Dirichlet} Let $G$ be a discrete subgroup of $M(\bar\R^n)$,
and $x$ a point in $\hyp^{n+1}$, which is not fixed by any
nontrivial element of $G$. The set $\{y \in
\hyp^{n+1}| d_h(y, x) \leq d_h(y, g(x))\; \forall g
\in G \}$ is called the Dirichlet domain for $G$ centred at $x$, where $d_h$
denotes the hyperbolic distance.
\end{defi}

It is easy to see that any Dirichlet domain is convex and the
interior of the intersection of the closure of a Dirichlet domain with $\bar
\R^n$ is a fundamental domain as defined before.

The following follows immediately from the definition of conical limit points.

\begin{lem}
\label{no conical limit point}
Let $D$ be a Dirichlet domain of  a discrete group $G \subset M(\bar\R^n)$.
Then $\bar D \cap \bar \R^n$ contains no conical limit points, where $\bar D$
denotes the closure of $D$ in $\B^{n+1}=\hyp^{n+1}\cup \bar \R^n$.
\end{lem}

Now, we consider how a Dirichlet domain of a geometrically finite
group intersects standard parabolic regions. We shall make use of
the following result of Bowditch \cite{Bow}. For a $G$-invariant set
$S$ on $\bar \R^n$, we say a collection of subsets $\{A_s\}_{s\in
S}$ is {\em strongly invariant} if $gA_s=A_{gs}$ for any $s \in S$
and $g \in G$, and $A_s \cap A_t = \emptyset$ for any $s\not=t \in
S$.
We should note that each $A_s$ is in particular precisely invariant
under $\mathrm{Stab}_G(s)$ in $G$ in the sense as defined before.

\begin{lem}
\label{disjoint horoballs} Let $\Pi$ be the set of all bounded
parabolic fixed points contained in the limit set $\Lambda(G)$  of a
discrete group $G \subset M(\bar\R^n)$. Then we can choose a
standard parabolic region $B_p$ at $p$ for each $p \in \Pi$ in
such a way that $\{B_p|p \in \Pi\}$ is strongly invariant.
\end{lem}

Using this lemma, we can show the following, which is essentially contained in
the argument of \S 4 in Bowditch \cite{Bow}.

\begin{pro}
\label{essentially finite} Let $D$ be a Dirichlet domain of a
geometrically finite group $G \subset M(\bar\R^n)$. Let $\{B_p\}$ be
the collection of standard parabolic regions obtained as in the
preceding lemma. Then there is a finite number of points $p_1, \dots
, p_k \in \bar D \cap \Pi$ such that $\bar D \setminus
\cup_{i=1}^k(\mathrm{Int} B_{p_i} \cup \{p_i\})$ is
compact and contains no limit point of $G$.
\end{pro}

\begin{proof}
Choose a family of standard parabolic regions $\{B_p\}$ as in
Lemma \ref{disjoint horoballs}. Since $G$ is geometrically
finite, every limit point of $G$ is either a conical limit point or
a  parabolic vertex. By Lemma \ref{no conical limit point}, no limit
point on $\bar D$ is a conical limit point. Therefore $\{B_p\}$
covers all limit points contained in $\bar D$.

Suppose that there are infinitely many distinct $B_{p_i}$ among
$\{B_p\}$ with $p_i \in \bar D$. By taking a subsequence, we can
assume that $\{p_i\}$ converges to a point $q \in \bar D$, which is
also contained in $\Lambda(G)$, hence in $\Pi$. By taking a
subsequence again, we can further assume that  all the $p_i$ belong
to either the same $\mathrm{Stab}_G(q)$-orbit or distinct
$\mathrm{Stab}_G(q)$-orbits. We first consider the former case. Let
$\alpha_i$ be the geodesic line connecting $p_i$ to $q$, which must
be contained in $D$. Since all $p_i$ belong to the same orbit, there
are $h_i \in \mathrm{Stab}_G(q)$ such that $h_i(p_i)=p_1$. By taking
a subsequence again, we can assume that all $h_i$ are distinct.
Then, the geodesic $\alpha_1$ is shared by infinitely many
translates of $h_iD$. This contradicts the local finiteness of
the tranlates of the Dirichlet domain $D$.

Since $q$ is a parabolic vertex, by Lemma \ref{bounded}, we see that
\linebreak$(\Lambda(G) \setminus \{q\})/\mathrm{Stab}_G(q)$ is compact.
Therefore, by taking a subsequence again, we can assume that there
are $g_i \in \mathrm{Stab}_G(q)$ such that $\{g_ip_i\}$
converges to a point $r \in \bar\R^n\setminus \{q\}$. We can
assume that all the $g_i$ are distinct by taking a subsequence. Let
$\alpha_i$ be the geodesic line connecting $p_i$ and $q$
as before. Then $g_i \alpha_i$ converges to the geodesic line
connecting $r$ to $q$. Since $g_i \alpha_i$ is contained in $g_i
D$, this again contradicts the local finiteness of the translates of $D$.
\end{proof}

Another easy consequence of Lemma \ref{disjoint horoballs} is the following.
\begin{cor}
\label{isometric sphere}
Let $G$ be a discrete subgroup of $M(\bar \R^n)$.
In the upper half-space model of $\hyp^{n+1}$, suppose that
$\infty$ is a parabolic vertex of $G$.
Then the Euclidean  radii of the isometric spheres $I(g)$ of $g\in
G-\mathrm{Stab}_G(\infty)$ are bounded from above.
\end{cor}
\begin{proof}
Consider the set of standard parabolic regions $\{B_p\}_{p \in \Pi}$
obtained by Lemma \ref{disjoint horoballs}. Since $\infty$ is a
bounded parabolic fixed point, a standard parabolic region
$B_\infty$ and its translates $g B_\infty$ by elements $g \in G -
\mathrm{Stab}_G(\infty)$ are among $\{B_p\}$. Let $B'_\infty$ be the
maximal  horoball contained in $B_\infty$. Then there is a
number $h$ such that $B'_\infty= \{(z_1, \dots ,
z_{n+1})|z_{n+1}\geq h\}\cup \{\infty\}$, which is equal to the
height of $\mathrm{Fr}B'_\infty$.

Fix an element $g \in G-\mathrm{Stab}_G(\infty)$. By enlarging
$B'_\infty$, we get a  horoball $B''$ which touches
$g^{-1} B''$ at one point. Let $h' < h$ be the height of
$\mathrm{Fr}B''$. Then the point $B'' \cap g^{-1}B''$ has height
$h'$. The isometric sphere $I(g)$ of $g$ must contain the point $B''
\cap g^{-1}B''$ since the reflection in $I(g)$ sends $g^{-1} B''$ to
$B''$. Therefore the Euclidean radius of $I(g)$ is equal to $h'$,
which is bounded above by the constant $h$ independent of $g$.
\end{proof}

This implies the following fact in the conformal ball model, which
is Corollary $G.8$ in Maskit \cite{Mas}.

\begin{cor}
\label{radii of isometric spheres} We regard $G$ as above as acting
on the ball $\B^{n+1}$ or $L=\bar \R^{n+1} \setminus \B^{n+1}$, and
let $p \in \partial \B^{n+1}=\partial L$ be a parabolic vertex of
$G$. Suppose that $g_n \in G$ are distinct elements. Then the radius
with respect to the ordinary Euclidean metric on $\B^{n+1}$ or $L$
of the isometric sphere $I(g_k)$ goes to $0$ as $k \rightarrow
\infty$.
\end{cor}

\section{Blocks}
Throughout this section, we assume that $G$ is a discrete subgroup of
$M(\bar \R^n)$ and $J$ is a subgroup of $G$.

\begin{defi}
A closed $J$-invariant set $B$, containing at lease two points, is
called a block, or more specifically $(J, G)$-block
 if it satisfies the following conditions.
\begin{enumerate}
\item  $B\cap \Omega(G)=B\cap \Omega(J)$, and $B\cap \Omega(G)$ is
precisely invariant under $J$ in $G$.

\item If $U$ is a peak domain for a parabolic fixed point $z$ of
$J$ with the rank of $\mathrm{Stab}_J(z)$ being $k<n$, then there is a
smaller peak domain $U'\subset U$ such that $U'\cap
Fr B=\emptyset$.
\end{enumerate}
\end{defi}

Let $S$ be a $(J, G)$-block, and let $S$ be a topological
$(n-1)$-dimensional sphere in $\bar \R^n$. Then $S$ separates $\bar
\R^n$ into two open sets.
We say that $S$ is {\it precisely embedded} in $G$ if  $g(S)$ is disjoint from one
of the two open sets for any $g \in G$.
\smallskip

A $(J, G)$-block is said to be {\it strong} if every parabolic fixed point
of $J$ is a parabolic vertex of $G$.
\medskip

Then we have the following.

\begin{thm}
\label{diameter vanishes} Suppose that $G$ is a discrete subgroup of
$M(\bar\R^n)$. Let $J$ be a geometrically finite subgroup of $G$ and
$B\subset\bar \R^n$  a $(J, G)$-block such that for every
parabolic fixed point $z$ of $J$ with the rank of $Stab_J(z)$ being
less than $n$, there is a peak domain $U_z$ for $J$ with $U_z\cap
B=\emptyset$. Let $G=\cup g_k J$ be a coset decomposition. Then we
have $\mathrm{diam}(g_k(B))\rightarrow 0,$ where $\mathrm{diam}(M)$
denotes the diameter of the set $M$ with respect to the ordinary
spherical metric on $\bar\R^n$.
\end{thm}

\begin{proof} By conjugating  $G$ by an element of $M(\bar \R^n)$, we can assume that
$\mathrm{Stab}_G(\mathbf{0})=\mathrm{Stab}_G(\infty)=\{id\}$ when we
regard $G$ as acting on $\bar \R^{n+1}$ by considering the
Poincar\'{e} extension. Let $L$ denote the exterior of $\B^{n+1}$
with the point $\infty$, which we regard also as a model of
hyperbolic $(n+1)$-space. Then $J$ is also geometrically finite as
a discrete group acting on  $L$. Let $P$ be a Dirichlet
domain for $J$ in $L$.

Let $g$ be some element of $G-J$. For a fixed $g$,  the
set $\{(g\circ j)^{-1}(\infty)=j^{-1}\circ g^{-1}(\infty)| j\in
J\}$  is $J$-invariant.
Then for each coset $g_k J$, we can choose a representative $g_k$ in such a way that
$a_k=g_k^{-1}(\infty)$, which is the centre of the isometric sphere of $g_k$, lies in $P$.

Now, by Proposition \ref{essentially finite}, there are finitely
many standard parabolic regions $B_{p_1}, \dots , B_{p_s}$ in $L$
around parabolic vertices $p_1, \dots , p_s$ on $\bar P$ such that
$\bar P \setminus \cup_i (\mathrm{Int} B_{p_i} \cup \{p_i\})$ is
compact and contains no limit point of $J$. We number them in such a
way that $\mathrm{Stab}^*_J(p_1), \dots , \mathrm{Stab}^*_J(p_r)$
have rank $n$ whereas $\mathrm{Stab}^*_J(p_{r+1}), \dots ,
\mathrm{Stab}^*_J(p_s)$ have rank less than $n$. We can assume that
for $j \geq r+1$, we have $B_{p_j} \cap \bar \R^n \cap B=\{p_j\}$
because of the following:
 By  our assumption in the theorem, we can make $B_{p_j}$ smaller so that it satisfies
 this condition.
 Also it is clear that for the old $B_{p_j}$, there is no limit point of $J$ in
 $\bar \R^n \cap B_{p_j}$
 other than $p_j$, which is also contained in the new $B_{p_j}$.
On the other hand no point in $\bar P$ can converge to
$p_j$ from outside this smaller $B_{p_j}$ since $p_j$ is not a
conical limit point, which implies that the compactness is
preserved.

For horoballs $B_{p_1}, \dots , B_{p_r}$, we have the following.
\begin{cl}
\label{disjoint}
We can choose the horoballs $B_{p_1},\cdots,B_{p_r}$
sufficiently small so that
$B_{p_i}\cap
G(\infty)=\emptyset$ for each $i$ ($1\leq i\leq r$).
\end{cl}
\begin{proof}
We identify $L$ with the standard upper half-space model of
hyperbolic $(n+1)$-space, which we denote by $\hyp^{n+1}$. By
conjugation, we can assume that  $e=(0, \dots ,0,1)$ corresponds to
$\infty \in L$ under the identification of $\hyp^{n+1}$ with $L$.
Regarding $G$ as acting on this $\hyp^{n+1}$ and $B_{p_1}, \dots,
B_{p_r}$ lying in $\B^{n+1}$, what we have to show is that
$B_{p_i}\cap G(e)=\emptyset$ for each $i$.

We shall show that how we can make $B_{p_1}$ satisfy this condition.
Conjugating $G$ by an isometry of $\hyp^{n+1}$, we may assume that
$p_1=\infty$. Then Corollary \ref{isometric sphere} implies
that the radii of the isometric spheres $I(g)$ of $g\in
G-\mathrm{Stab}_G(\infty)$ are bounded from above by some constant
$r_0$. We set  $B_{p_1}=\{x\in \hyp^{n+1}| x_{n+1} \geq
2\max\{1, r_0^2\}\} \cup \{\infty\}.$

Any $h\in \mathrm{Stab}_G(\infty)$ can be represented as a
transformation of $\R^n$ in the form  $h(x)=Ax+b$ for $A\in O(n)$
and $b\in \R^n$. Let $\tilde{h}$ denote $h$ regarded as an isometry
of $\hyp^{n+1}$. Then we have $\tilde{h}(e)=(b,1)$, hence
$\tilde{h}(e)\notin B_{p_1}$.

For any $\displaystyle g\in G-\mathrm{Stab}_G(\infty)$,
let $r_g$ denote the radius of the isometric sphere $I(g)$.
Then $g(x)$ is represented as a transformation of $\bar \R^n$ in the form
$a+\frac{r_g^2A(x-b)}{|x-b|^2}$ for some $A\in O(n)$ and $a, b\in
\R^n$ (see \cite{Aga} or \cite{Bea}).
As before we denote by $\tilde{g}$ the transformation $g$ regarded as an
isometry of $\hyp^{n+1}$.
Then we have
$$\tilde{g}(e)=(a-\frac{r_g^2Ab}{|b|^2+1},\frac{r_g^2}{|b|^2+1})
\;\; \mbox{and}\;\; \frac{r_g^2}{|b|^2+1}\leq r_0^2,$$ which implies
that  $\tilde{g}(e)\notin B_{p_1}$. We make each $B_{p_i}$ smaller
in the same way. It is clear that even after changing the horoballs,
$\bar P \setminus \cup_i (\mathrm{Int} B_{p_i} \cup
\{p_i\})$ is compact and contains no limit point of $J$
since $B_{p_j}$ intersects $\bar P \cap \bar \R^n$ only at $p_j$
$(1\leq j\leq r)$ and $p_i$ is not a conical limit point.
\end{proof}

%

Recall that  $a_k = g_k^{-1}(\infty)$ is in $P$. By taking
a subsequence, we have only to consider the cases when every $a_k$
lies outside all the standard parabolic regions $B_{p_j}$ and when
all the $a_k$ lie in some $B_{p_j}$.

First consider the case when  every $a_k$ lies outside the
$B_{p_j}$. Since $a_k\in \bar{P}$ and $ \bar P\setminus
\cup (\mathrm{Int}(B_{p_j}) \cup \{p_j\})$ is
compact, the sequence $\{a_k\}$ converges to a point $x\in
\bar P\setminus \cup (\mathrm{Int}(B_{p_j})\cup
\{p_j\})$. Suppose that $x$ is contained in $B$. Then $x$
must lie in $B\cap\Lambda(G)=B\cap\Lambda(J)$, which contradicts the
fact that $ \bar P\setminus \cup (\mathrm{Int}(B_{p_j}) \cup
\{p_j\})$ contains no limit point of $J$. Therefore, it follows that
the $a_k$ are uniformly bounded away from $B$. Since the $g_k$ are
distinct elements,  the radius with respect to the Euclidean metric
of the conformal ball model of the isometric sphere $I(g_k)$
converges to $0$ by Corollary \ref{radii of isometric spheres}.
Therefore, we see that  $B$ lies outside  the isometric sphere
$I(g_k)$ for sufficiently large $k$. This means $g_k(B)$ lies inside
the isometric sphere $I(g_k^{-1})$. This implies that
$\mathrm{diam}(g_k(B))\rightarrow 0.$

Next we consider the case when the $a_k$ lie in some standard
parabolic region $B_{p_j}$. By Claim \ref{disjoint}, we see that
$B_{p_j}$ is not a horoball; hence $B_{p_j}$ is an extended
horoball, i.e., $j\geq r+1$. Furthermore, if  $\{a_k\}$ does not
converge to $p_j$, then  we can take $B_{p_j}$ smaller. Therefore,
we can assume that $\{a_k\}$ converges to $p_j$.

By composing a rotation of the sphere $\bar \R^n$, we may assume that
$p_j$ is at the north pole $(0, \dots ,0,1)$.
Let $S$ be the $n$-sphere of radius $1$ centred at $p_j$, and let
$\phi$ be the reflection in $S$.
Let
$B^{'}\subset B_{p_j}$ be the largest horoball contained in $B_{p_j}$
touching $\bar \R^n$ at $p_j$.

We denote points in $\R^{n+1}$ as $(\mathbf{z},t)$ with
$\mathbf{z}\in \R^n$ and $t\in \R$. Then we have
$p_j=(\mathbf{0},1)$.Take $B_{p_j}$ to be small enough so that
$B'=\{(\mathbf{z},t)||\mathbf{z}|^2+(t-s'-1)^2\leq  s'^2\}$
for some $s'$ satisfying $0<s'<\frac{1}{2}$, and
$$\phi(z,t)=(\frac{\mathbf{z}}{|\mathbf{z}|^2+(t-1)^2},\;\;
\frac{|\mathbf{z}|^2+t^2-t}{|\mathbf{z}|^2+(t-1)^2}).$$

We deduce that
$$\phi(\B^{n+1})=\{(\mathbf{z},t)|
t\leq\frac{1}{2}\}\cup\{\infty\}\;\;\mbox{and}\;\;
\phi(B')=\{(\mathbf{z},t)| t\geq
1+\frac{1}{2s'}\}\cup\{\infty\}.$$

For any $j\in \mathrm{Stab}_J(p_j)$, we have $\phi j\phi(\infty)=\infty$
since $\phi(\infty)=p_j$.
Consider the decomposition $\R^{n+1}= \R^m \times \R^{n-m}\times \R$,
where $m (<n)$ is the
rank of $\mathrm{Stab}_J(p_j)$.
Let $\phi j\phi(z)= U(z)+\mathbf{a} $ be an
arbitrary element of $\phi \mathrm{Stab}_J(p_j)\phi$,
where $U$ denotes a rotation.
By Theorem \ref{Bieberbach}, we may
assume that the rotation $U$ leaves $\R^m$ and $\R^{n-m}$ invariant
and the vector $\mathbf{a}$ lies in the subspace $\R^m$.
Also, if $\phi j\phi\in
\phi \mathrm{Stab}^*_J(p_j)\phi$, then its restriction to
the subspace $\R^m$ is
a translation.
Hence, we have
$$\phi(B_{p_j})=\{(\mathbf{z},t)| \sum_{i=m+1}^n z_i^2+t^2\geq (1+\frac{1}{2s'})^2,
t\geq \frac{1}{2}\}\cup\{\infty\},$$ where $z_i$ denotes the $i$-th component of
$\mathbf{z}$.

Since $B_{p_j}\cap B=\{p_j\}$, we have \begin{equation}
\label{B} \phi(B)\subset \{(\mathbf{z},t): \sum_{i=m+1}^n
z_i^2+\frac{1}{4}< (1+\frac{1}{2s'})^2, t=\frac{1}{2}\}\cup\{\infty\}.
\end{equation}
We should  recall that $\phi \mathrm{Stab}^*_J(p_j) \phi$ acts on $\R^m$ cocompactly.
Therefore, we can take representatives $g_k$ so that the projections of
$\phi(a_k)=\phi(g_k^{-1}(\infty))$ to $\R^m$ stay within a compact subset
of $\R^m$ by multiplying elements of $\mathrm{Stab}^*_J(p_j)$ to the original $g_k$.
Note that by changing representatives, we do not have the condition that
$a_k \in  P$ any more, but still the $a_k$ are contained in $B_{p_j}$.
This means that there is a constant $L$ such that
$\phi(a_k)\in \{(z,t)| \sum_{i=1}^m z_i^2 < L,
t>\frac{1}{2}\}\cap\phi(B_{p_j}).$

\begin{cl}
\label{Lipschitz} There is a constant $K>0$ such that for every
$a_k\in B_{p_j}$ and every $y\in B$, we have $|a_k-y|\geq
K|a_k-p_j|$.
\end{cl}
\begin{proof}
Suppose, seeking a contradiction, that such a $K$
does not exist.
Then there exist a sequence $\{y_s\}\subset B$ and
a subsequence $\{a_{k_s}\}$ of $\{a_k\}$ such that
\begin{equation}
\label{go to 0}
\frac{|a_{k_s}-y_s|}{|a_{k_s}-p_j|}\to 0\;\;\mbox{as}\;\;s\to \infty.
\end{equation}
We shall denote $a_{k_s}$ by $a_s$ for simplicity.

We can assume that $y_s\not= p_j$ for all $s$.
Then, since
$$|\phi(a_s)-\phi(y_s)|=\frac{|a_s-y_s|}{|y_s-p_j||a_s-p_j|}\;\;\mbox{and}\;\;
|\phi(y_s)-p_j||y_s-p_j|=1,$$ we have

\begin{equation}
\label{sum}
\begin{split}
\frac{|a_s-y_s|^2}{|a_s-p_j|^2}&=\frac{|\phi(a_s)-\phi(y_s)|^2}{|\phi(y_s)-p_j|^2}\\
&=\frac{\sum^{m}_{i=1}
(\phi(a_s)-\phi(y_s))^2_i+\sum^{n+1}_{i=m+1}(\phi(a_s)-\phi(y_s))^2_i}
{\sum^m_{i=1}(\phi(y_s))^2_i+\sum^{n+1}_{i=m+1}
(\phi(y_s)-p_j)^2_i}.
\end{split}
\end{equation}

We shall show  that there exists $M>0$ such that

\begin{enumerate}
\item $\sum^{m}_{i=1}(\phi(a_s))^2_i\leq M$ for all $s$;
\item $\sum^{n+1}_{i=m+1}(\phi(y_s)-p_j)^2_i\leq M$ for all
$s$; and
\item
$\sum^{n+1}_{i=m+1}(\phi(a_s)-\phi(y_s))^2_i\to\infty$ as $s\to
\infty$.
\end{enumerate}
The inequality (1) follows from the fact that we choose $a_k$ so
that the projections of $\phi(a_k)$ to $\R^m$ stay in a compact
subset. The second one is a consequence of (\ref{B}). We now turn to
the third inequality. Since $\{a_s\}$ was assumed to converge to
$p_j$, we see that $\phi(a_s)$ tends to $\infty$, which means that
$\sum_{i=1}^{n+1}(\phi(a_s))_i^2 \to \infty$. On the other hand,  we
know that $\sum_{i=1}^m(\phi(a_s))^2_i \leq M$ by (1), and that
$\sum_{i=m+1}^{n+1}(\phi(y_s))_i^2$ is bounded above independently
of $s$ by (2). These imply  (3).

Then (\ref{go to 0}), (\ref{sum}), (2) and (3) imply that
$$\sum^{m}_{i=1}(\phi(y_s))^2_i\to \infty\;\;\mbox{as}\;\;s\to\infty.$$
It follows from (1) that for all sufficiently large $s$,
$$\frac{|a_s-y_s|}{|a_s-p_j|}\geq \frac{1}{2}.$$
This is a contradiction and we have completed the proof of Claim \ref{Lipschitz}.
\end{proof}

Let $\rho_k$ be the Euclidean radius of the isometric sphere of $g_k$ in $L$. Then
we have the following.
\begin{cl}
\label{radius}
 If all  $a_k$ lie inside the extended
horoball $B_{p_j}$, then we have
$\displaystyle\frac{{\rho^2_k}}{|a_k-p_j|}\rightarrow 0$.
\end{cl}

\begin{proof}
Suppose that there is  $\delta>0$ such that
$\displaystyle\frac{{\rho^2_k}}{|a_k-p_j|}\geq\delta$.
Then
$\displaystyle|g_k(p_j)-g_k(\infty)|=
\frac{{\rho^2_k}}{|a_k-p_j|}\geq\delta.$

We can apply Proposition \ref{conical limit} by identifying $L$ with
$\B^{n+1}$ by the reflection in $\partial \B^{n+1}$ and taking into
account the fact that the Euclidean metric does not distort much by
the reflection  near $\partial \B^{n+1}$\, and  see that
$p_j$  is a conical limit point of $G$. This contradicts
Lemma \ref{no conical limit point} since $p_j$ lies in $\bar P$.
\end{proof}

We shall conclude the proof of Theorem \ref{diameter vanishes}. Let
$\delta_k$ be the distance from $a_k$ to $B$. Since $\delta_k$ is
the infimum of $|a_k-y|$ for $y\in B$,  by Claim \ref{Lipschitz}, we
have $\delta_k\geq K|a_k-p_j|$. Since Proposition I.C.7 in
\cite{Mas4} holds for $g\in M(\bar R^n)$, we have
$$\mathrm{diam}(g_k(B))\leq \frac{2\rho^2_k}{\delta_k}\leq
\frac{2K^{-1}\rho^2_k}{|a_k-p_j|}.$$
This implies that $\mathrm{diam}(g_k(B))\to 0$ by Claim \ref{radius}.
\end{proof}

\section{The  Combination Theorem }
%
%
%
In this section, we shall state and prove our main theorem, which is a combination
theorem for discrete groups in $M(\bar \R^n)$.
Before that we shall prove the following lemma which constitutes the key step for
the proof of our main theorem.

\begin{lem}
\label{key lemma} Let $G_1$ and $G_2$ be discrete subgroups of
$M(\bar \R^n)$. Suppose that $J$ is a subgroup of $G_1 \cap G_2$,
which coincides with neither $G_1$ nor $G_2$. Suppose that there is
a topological $(n-1)$-sphere $S$ dividing $\bar \R^n$ into two
closed balls $B_1$ and $B_2$ such that each $B_m$ is a
$(J,G_m)$-block. Let $D_1, D_2$ be  fundamental
sets for $G_1$ and $G_2$, respectively such that
$J(D_m\cap B_m)=B_m\cap$$^\circ$$\Omega(J)$ for $m=1,2$, and
$D_1\cap S=D_2\cap S$. Set $D=(D_1\cap B_2)\cup(D_2\cap B_1)$ and
$G=<G_1, G_2>$. Then the following hold.

\begin{enumerate}
\item $S$ is also a $(J,G_m)$-block for $m=1,2$.

\item
$S\cap\Lambda(G_1)=S\cap\Lambda(G_2)=S\cap\Lambda(J)=\Lambda(J)$.

\item Both $G_1$ and $G_2$ have non-empty regions of discontinuity, and $B_m^\circ$
is contained in $\Omega(G_m)$ for $m=1,2$, where $B^\circ_m$ is the
interior of $B_m$ in $\bar \R^n$.

\item $B^\circ_m$ is precisely invariant  under $J$ in $G_m$.

\item For any $g\in G_m-J$ $(m=1,2)$, $g(B_m)\cap B_m=g(S)\cap
S\subset \Lambda(G_m)$.

\item For any $g\in G_m$, we have $g(D_m\cap B_{3-m})\subset B_{3-m}$ and
$g(D_m\cap B^\circ_{3-m})\subset B^\circ_{3-m}$.

\item Let $G_m=\bigcup g_{km}J$ be a coset decomposition for $m=1,2$. If $J$
is geometrically finite, then $\mathrm{diam}(g_{km}(B_m))\rightarrow
0$ as $k \rightarrow \infty$ where $\mathrm{diam}$ denotes
the diameter with respect to the ordinary spherical metric on $\bar
\R^n$.

\item $(B^\circ_1, B^\circ_2)$ is an interactive pair.

\item If $\Lambda(J)\neq\Lambda(G_1)$ or
$\Lambda(J)\neq\Lambda(G_2)$, then $(B^\circ_1, B^\circ_2)$ is a
proper interactive pair.

\item  If $D\neq\emptyset$ and $J$ is geometrically finite, then
$(B^\circ_1, B^\circ_2)$ is a proper
interactive pair.
\end{enumerate}
\end{lem}

\begin{proof}
(1). This is obvious since $S$ is contained in $B_m$.

(2). By Lemma \ref{lem:limit set}, we see that  $\Lambda(J)$ is contained
in $S$; hence $S \cap \Lambda(J)=\Lambda(J)$.
Since $S$ is a $(J, G_m)$-block for $m=1,2$ by (1), we have
$S\cap\Lambda(G_m)=S\cap\Lambda(J)$.
This implies (2).

(3). Since $\Lambda(J)$ is contained in  $S$, we see that $B^\circ_m
\cap \Omega(J)= B_m^\circ$. On the other hand, since $B_m$
is a $(J,G_m)$-block, we have $B^\circ_m\cap
\Omega(G_m)=B^\circ_m\cap \Omega(J)=B^\circ_m\neq \emptyset$. Thus
both $G_1$ and $G_2$ have non-empty regions of discontinuity and
$\Omega(G_m)$ contains $B_m^\circ$.

(4).
 Since
 $B^\circ_m\subset\Omega(G_m)$, by the definition of  blocks,
 $B_m\cap\Omega(G_m)$ is precisely invariant under $J$ in $G_m$, and
$J(S)=S$, we see that $B_m^\circ$ is precisely invariant under $J$ in $G_m$.

(5).
Since $B_m \cap \Omega(G_m)$ is precisely invariant under $J$
in $G_m$, for every  $g\in G_m-J$,
$g(B_m\cap \Omega(G_m))\cap (B_m\cap \Omega(G_m))=\emptyset$.
It follows
$(g(B_m)\cap B_m)\cap\Omega(G_m)=\emptyset.$
Then we see that  (4) implies (5).

(6).
For any $j\in J\subset G_m$,
$j(D_m\cap B_{3-m})\subset
j(B_{3-m})=B_{3-m}\;\ \mbox{and}\;\ j(D_m\cap B^\circ_{3-m})\subset
j(B^\circ_{3-m})=B^\circ_{3-m}$.
Hence we have only to consider the case when $g$ lies in $G_m-J$.
Suppose that there exists an element $g\in G_m-J$ such that
$g(D_m\cap B_{3-m})\cap B_{m}\neq\emptyset$.
Take points  $x\in g(D_m\cap B_{3-m})\cap B_{m}$ and
$y\in D_m\cap B_{3-m}$ such that $x=g(y).$
Since $x$ lies in $B_m \cap g(D_m \cap B_{3-m})  \subset B_{m}\cap$
$^\circ$$\Omega(G_m)\subset B_m\cap$$^\circ$$ \Omega(J)=J(D_m\cap B_m)$,
there are an element $j\in J$ and a point $z\in D_m\cap B_{m}$ such that $j(z)=x$.
Then
$j(z)=g(y)$. Since $z$ and $y$ are $G_m$-equivalent points of $D_m$,
we have $z=y$ and $j=g$, which is a contradiction. Therefore, for any $g\in
G_m-J$, we have $g(D_m\cap B_{3-m})\cap B_{m}=\emptyset$ and $g(D_m\cap
B_{3-m})\subset B^\circ_{3-m}$.
Thus we have proved (6).

(7). By (1),  we know that  $S$ is a $(J, G_m)$-block. Also we
should note that since $\mathrm{Fr}S=S$, by the definition of
blocks, for any parabolic vertex $z$ of $J$ on $S$ with the
rank of $Stab_J(z)$ being less than $n$, there is a peak
domain centred at $z$ which is disjoint from $S$, and that every
parabolic fixed point is a parabolic vertex if $J$ is geometrically
finite. Therefore by Theorem \ref{diameter vanishes},
$\mathrm{diam}(g_{km}(S))\rightarrow 0$ as $k \rightarrow \infty$.
On the other hand since $B_m$ is a $(J,G_m)$-block,
$\mathrm{diam}(g_{km}(S)) \rightarrow 0$ implies
$\mathrm{diam}(g_{km}(B_m))\to 0$, and we have completed the proof
of (7).

(8).
This follows from (4) and Proposition \ref{Maskit}.

(9).
 If $(B^\circ_1, B^\circ_2)$ is not
proper, then $B^\circ_1\cup
B^\circ_2=G_1(B^\circ_1)\subset\Omega(G_1)$ and $B^\circ_1\cup
B^\circ_2=G_2(B^\circ_2)\subset\Omega(G_2)$. It follows that for
each  $m$, we have $\Lambda(G_m)\subset S$.
On the other hand, by (2),
we have $\Lambda(G_m)=S\cap\Lambda(G_m)=S\cap\Lambda(J)=\Lambda(J)$.
Therefore if one of $\Lambda(G_1), \Lambda(G_2)$ is not equal to $\Lambda(J)$,
then $(B_1^\circ, B_2^\circ)$ is a proper interactive pair.

(10). Suppose that $D$ is non-empty and $J$ is geometrically finite.
Then we can  assume that $D_1\cap B_2\neq\emptyset$, for the case
$D_2 \cap B_1$ can be proved just by interchanging the indices. We
divide the argument into two cases: the case when $D_1 \cap S \neq
\emptyset$ and the one when $D_1 \cap B_2^\circ \neq \emptyset$.

Suppose first that there is a point $x\in D_1\cap
S=D_2\cap S$.
Recall that $D_1$ is contained in $\Omega(G_1)$, and that for $g\in G_1-J$,
we have $g(B_1) \cap B_1 \subset \Lambda(G_1)$ by (5).
These imply that no $(G_1-J)$-translates of $B_1$ pass
through $x \in D_1 \cap S \subset D_1 \cap B_1$.
By the same argument, we see that no
$(G_2-J)$-translates of $B_2$ pass through $x$.

Next we shall show that $(G_m-J)(B_m)$ cannot accumulate at $x$.
First we should note that the translate of $B_m$ by an element  of
$G_m$ depends only on the cosets of $G_m$ over $J$ since
$J$ stabilises $B_m$. Suppose that $(G_m-J)(B_m)$
accumulates at $x$. Then there are elements $g_k$ in
$G_m-J$, which we can assume to belong to distinct cosets,
and points $y_k \in B_m$ such that $\{g_k(y_k)\}$ converges to $x$.
Since we assumed that $J$ is geometrically finite, by (7) we see
that $\mathrm{diam}(g_k(B_m)) \rightarrow 0$. Therefore if we choose
one point $y$ in $B_m$, then $\{g_k(y)\}$ also converges to $x$.
 This means that $x$ is a limit point of $G_m$, which contradicts the
 assumption that $x$ lies in $D_m$.

By these two facts which we have just proved, we see that there is a
neighborhood of $x$ which is disjoint from
$(G_m-J)(B_m)$ for each $m$.
This implies in particular that there is a point in $B_{3-m}^\circ$ which
is not contained in the $G_m$-translates of $B_m$.
Hence, in this case, $(B^\circ_1,
B^\circ_2)$ is proper.

Now we assume that there is a point $x\in D_1\cap B^\circ_2$.
If $x\in (G_1-J)(B^\circ_1)$, then there are an
element $g\in G_1-J$ and a point $y\in B^\circ_1$ with $x=g(y)$.
Since $y$ lies in $B^\circ_1\cap
$$^\circ\Omega(G_1)\subset B^\circ_1\cap
$$^\circ\Omega(J)=J(D_1\cap B^\circ_1)$, there are
an element $j\in J$ and  a point $z\in D_1\cap B^\circ_1$
with $y=j(z)$, which implies $x=gj(z)$.
Since $D_1$ is a fundamental set of $G_1$, it follows that $x=z$ and $g=j^{-1}$,
which is a contradiction. Therefore $x$ is not contained in $(G_1-J)(B^\circ_1)$ and
$(B^\circ_1, B^\circ_2)$ is proper.
Thus we have proved (10).
\end{proof}

\begin{defi}
Let $\{S_j\}$ be a collection of topological
$(n-1)$-spheres. We say that the sequence $\{S_j\}$ {\it nests
about} the point $x$ if the following are satisfied.
\begin{enumerate}
\item The spheres $S_j$ are pairwise disjoint.

\item For each $j$, the sphere $S_j$ separates
$x$ from the precedent $S_{j-1}$;

\item For any point $z_j\in S_j$, the sequence $\{z_j\}$ converges to $x$.
\end{enumerate}
\end{defi}

Now we can state and prove our main theorem.

\begin{thm}
\label{main}
Let $J$ be a geometrically finite proper subgroup of two discrete
groups $G_1$ and $G_2$ in $M(\bar \R^n)$.
Assume  that there is a topological $(n-1)$-sphere $S$ dividing $\bar \R^n$ into
two closed balls $B_1$ and $B_2$ such that each $B_m$ is a $(J, G_m)$-block
and $(B^\circ_1, B^\circ_2)$ is a proper interactive pair.
Let $D_m$ be a fundamental set for $G_m$ for $m=1,2$  such that
$J(D_m\cap B_m)=B_m\cap$$^\circ\Omega(J)$,
$D_m\cap B_{3-m}$ is either empty or has nonempty interior, and
$D_1\cap S=D_2\cap S$.
Set $D=(D_1\cap B_2)\cup(D_2\cap B_1)$ and
$G=<G_1, G_2>$.
Then the following hold.

\begin{enumerate}
\item $G=G_1*_JG_2$.
\item $G$ is discrete.
\item If an element $g$ of $G$ is not loxodromic, then one of the following must hold.
\begin{enumerate}
\item $g $ is  conjugate to an element
of either $G_1$ or $G_2$.
\item $g$ is parabolic and is conjugate to an
element fixing a parabolic fixed point of $J$.
\end{enumerate}
\item  $S$ is a precisely embedded $(J, G)$-block.
\item  If $\{S_k\}$ is a sequence of distinct $G$-translates of
$S$, then $\mathrm{diam}(S_k)\rightarrow 0$, where $\mathrm{diam}$ denotes
the diameter with respect to the ordinary spherical metric on $\bar \R^n$.
\item There is a sequence of distinct $G$-translates of $S$
nesting about the point $x$ if and only if $x$ is a limit point of
$G$ which is not $G$-equivalent to a limit point of either $G_1$
or $G_2$.
\item $D$ is a fundamental set for $G$.
If both $D_1$ and $D_2$ are constrained, and $S\cap\mathrm{Fr} D$
consists of finitely many
connected components the sum of whose $(n-1)$-dimensional measures on $S$
vanishes, then $D$ is also constrained.
\item  Let $Q_m$ be the union of the $G_m$-translates of $B^\circ_m$, and
let $R_m$ be the complement of $Q_m$ in $\bar \R^n$.
Then $\Omega(G)/G=(R_1\cap\Omega(G_1))/G_1\cup(R_2\cap\Omega(G_2))/G_2,$
where the latter  two possibly disconnected orbifolds are identified along
their common possibly disconnected or empty boundary
$(S\cap \Omega(J))/ J$.
\end{enumerate}
\medskip
Furthermore, if each $B_m$ is precisely invariant under $J$ in $G_m$ for
$m=1, 2$, then the following hold.
\medskip
\begin{enumerate}
\setcounter{enumi}{8}
\item  $S$ is a strong $(J,G)$-block if and only if each $B_m$ is
a strong $(J, G_m)$-block.
\item If both $B_1$ and $B_2$ are strong, then, except for
$G$-translates of limit points of $G_1$ or $G_2$, every limit point of
$G$ is a conical limit point.
\item  $G$ is geometrically finite if and only if both $G_1$ and $G_2$
are geometrically finite.
\end{enumerate}
\end{thm}

\begin{proof}[Proof of (1)] Since $(B^\circ_1, B^\circ_2)$ is
proper, (1) follows from Theorem \ref{free product}.
\end{proof}

\begin{proof}[Proof of (2)]
Suppose that $G$ is not discrete. Then there is a sequence $\{g_k\}$
of distinct elements of $G$ which converges to the
identity uniformly on compact subsets. Express $g_k$ in a normal
form $g_k=\gamma_{n_k}\circ \gamma_{n_{k-1}}\circ\cdots \circ
\gamma_{n_1}$. We may assume that each $g_k$ has even length, for if
$g_k$ has odd length, then by Lemma \ref{interactive}, either
$g_k(B^\circ_1)\subset B^\circ_2$, or $g_k(B^\circ_2)\subset
B^\circ_1$, and such elements cannot converge to the identity. By
interchanging $B_1$ and $B_2$ if necessary, we may assume that
$(G_1-J)(B^\circ_1)$ is a proper subset of $B^\circ_2$ since
$(B_1^\circ,B_2^\circ)$ is proper. By choosing a subsequence, we may
assume that all the $g_k$ are $(1, 2)$-forms or all of them are $(2,
1)$-forms. It suffices to prove the case that every $g_k$ is a
$(1,2)$-form since if $g_k$ is a $(2,1)$-form, then $g_k^{-1}$ is a
$(1,2)$-form.

Since we assumed that each $g_k$ is a $(1, 2)$-form, we have
$g_k(B^\circ_2)\subset \gamma_{n_k}\circ
\gamma_{n_{k-1}}(B^\circ_2).$ If
$\gamma_{n_{k-1}}(B^\circ_2)=B^\circ_1$, then $g_k(B^\circ_2)\subset
\gamma_{n_k}(B^\circ_1)\subset B^\circ_2,$ with the last
inclusion being proper, and  if
$\gamma_{n_{k-1}}(B^\circ_2)$ is a proper subset of
$B^\circ_1$, then $g_k(B^\circ_2)\subset \gamma_{n_k}\circ
\gamma_{n_{k-1}}(B^\circ_2) \subset
\gamma_{n_k}(B^\circ_1)\subset B^\circ_2,$ with the last two
inclusions being proper. Therefore, in either case, we have
$g_k(B^\circ_2)\subset \gamma_{n_k}(B^\circ_1)\subset
B^\circ_2,$ with the last inclusion being proper. Thus
$B^\circ_2-g_k(B^\circ_2)\supset B^\circ_2-
\gamma_{n_k}(B^\circ_1)\supset
B^\circ_2-(G_1-J)(B^\circ_1)$. Since $g_k\rightarrow id$ on $B_2$
and $B^\circ_2\setminus (G_1-J)(B^\circ_1)\neq \emptyset$, this is a
contradiction.
\end{proof}

Now for a normal form $g=g_n\cdots g_1\in G$, we call $g$ {\it
positive} if $g_1\in G_1-J$ and we express it as $g>0$; we call $g$
{\it negative} if $g_1\in G_2-J$ and we express it as $g<0$.

Using this distinction, we consider a coset decomposition of $G$:
$$G=J\cup(\cup_{n,k} a_{nk}J)\cup(\cup_{n,k} b_{nk}J),$$
where $|a_{nk}|=|b_{nk}|=n,\ a_{nk}>0,\ \mbox{and}\ b_{nk}<0.$
Following
Apanasov \cite{Apa}, we set
$T_n=(\cup_k a_{nk}(B_1))\cup(\cup_k b_{nk}(B_2)),\ C_n=\bar \R^n\setminus T_n,\
C=\cup C_n,\ \mbox{and}\ T=\bar \R^n \setminus C=\cap T_n.$

Then we have the following.
\begin{lem}
\label{decreasing T}
$\{T_n\}$ is a decreasing sequence with respect to the inclusion,
that is, $T_1 \supset T_2 \supset \dots  $.
\end{lem}

\begin{proof}
Take a point $x\in T_n$ $(n>1)$. Then either there are an
element $a_{nk}>0$ with length $n$ and a point $y\in B_1$ satisfying
that $x=a_{nk}(y)$, or there are an element $b_{nk}<0$ with length
$n$ and a point $y\in B_2$ satisfying that $x=b_{nk}(y)$. In the
former case, if we express $a_{nk}$ in a normal form as
$g_n\circ\cdots\circ g_1$, then $g_1\in G_1-J$. Since $g_1(y)$ lies
in  $g_1(B_1)\subset B_2$, there is a point $z\in B_2$ with
$g_1(y)=z$. Therefore,  $x=a_{nk}(y)=g_n\circ\cdots\circ g_2(z)\in
b_{(n-1)s}(B_2)\subset T_{n-1}.$ In the latter case, by the same
argument we have $x\in T_{n-1}.$
\end{proof}

\begin{lem}
\label{precisely embedded}
The sphere $S$ is precisely embedded in $G$. If $S$ is precisely invariant
under $J$ in $G_1$ and $G_2$, respectively, then $S$ is precisely
invariant under $J$ in $G$.
\end{lem}

\begin{proof}
We shall first show that $S$ is precisely embedded.
 For any $g\in G$ with $|g|=0$, we have $g(S)=S$ and is disjoint
 from both $B_1^\circ$ and $B_2^\circ$.
If $|g|=1$, then $g\in G_m-J$ $(m=1,2)$, and $g(S)=g(\mathrm{Fr}
B_m)\subset g(B_m)\subset B_{3-m}.$
This means that $g(S)$ is disjoint from $B_m^\circ$.

Now let $g=g_n\circ\cdots\circ g_1$ be an $(m,k)$-form with $|g|>1$.
Then $g(S)=g(\mathrm{Fr} B_k)\subset g(B_k)\subset B_{3-m}$
since $g(B^\circ_k)\subset B^\circ_{3-m}$ by Lemma \ref{interactive}.
This means that $g(S)$ is disjoint from $B_m^\circ$ again, and
we have thus shown that $S$ is precisely embedded in $G$.

Now suppose that  $S$ is precisely invariant under $J$ both in $G_1$
and $G_2$. Since, as was shown above, for $g \in J$, we have
$g(S)=S$, we have only to show that $g(S) \cap S =\emptyset$ for $g
\in G-J$. Note that $g(S)=g(\mathrm{Fr} B_m)\subset g(B_m)\subset
B^\circ_{3-m}$ for any $g\in G_m-J$. Therefore, it remains to
consider the case when $|g|>1$. If $g=g_n\circ\cdots\circ g_1$ is an
$(m,k)$-form with $|g|> 1$, then $h=g^{-1}_n\circ g$ is a
$(3-m,k)$-form. It follows from Lemma \ref{interactive} that
$g(S)=g_n\circ h(S)=g_n\circ h(\mathrm{Fr} B_k)\subset g_n\circ h
(B_k)\subset g_n(B_m)\subset B^\circ_{3-m}$. Thus, we have shown
that for any $g\in G-J$, $g(S)\cap S=\emptyset.$
\end{proof}

\begin{lem}\label{in C_1}
$D\subset C_1$.
\end{lem}

\begin{proof}
We assume that $D\neq\emptyset$.
By interchanging $B_1$ and $B_2$ if necessary, we can assume that
$D_1\cap B_2\neq\emptyset$.
If there is a point $x\in D_1\cap S=D_2\cap S$, then  no
$(G_m-J)$-translates of $B_m$ pass through $x$ as was shown in the proof of
Lemma \ref{key lemma}-(10). This implies
that $x\in C_1$.

It remains to consider the case when $x\in D_1\cap B^\circ_2$.
 If $x\in (G_1-J)(B_1)$, then there are an element $g\in G_1-J$ and a point $y\in B_1$ with $x=g(y)$.
Since $y\in$$^\circ\Omega(G_1)\cap B_1\subset$$^\circ\Omega(J)\cap
B_1$, there are an element $j\in J$ and a point $z\in D_1\cap B_1$
with $y=j(z)$ by the assumption that $J(D_1 \cap
B_1)=$$^\circ$$\Omega(J) \cap B_1$ in Theorem \ref{main}.
Therefore we have  $x=gj(z)$, which implies that $x=z$ and $gj=id$.
This contradicts the assumption that $g$ lies in $G_1-J$.
Thus we have shown that $x\in C_1$.
\end{proof}

\begin{lem}
\label{D}
$D$ is contained in $^\circ\Omega(G),$ and precisely invariant under $\{id\}$ in $G$.
\end{lem}

\begin{proof}
 We shall  first prove that  $D$ is contained in  $\Omega(G)$.
 Suppose, on the contrary, that there is a point $z$ in $D\cap\Lambda(G)$.
 Since $D=(D_1\cap B_2)\cup (D_2\cap B_1)$, we can assume that $z\in D_1\cap
B_2$ by interchanging the indices if necessary.

\begin{cl}
\label{z}
 In this situation, we have $z\in D_1\cap S$.
\end{cl}
\begin{proof}[Proof of Claim \ref{z}]
Suppose not.
Then $z$ must be contained in  $D_1\cap B^\circ_2$.
 Since $z\in \Lambda(G)$,
it follows from Lemma \ref{source-sink} that there is a sequence
$\{g_k\}$ of distinct elements in $G$ such that  $g_k(y)\rightarrow
z$ for all $y$ with at most one exception. Since $z\in
B^\circ_2\subset \Omega(G_2)$ (by Lemma \ref{key lemma}-(3)) and
$z\in D_1 \subset\Omega(G_1)$, we have $|g_k|>1$, and we can assume
that each $g_k$ is a $1$-form. Since $g_k(B)\subset T_1$ for  $B$
which is equal to $B_1$ or $B_2$, Lemma \ref{in C_1}
implies that $z\in\mathrm{Fr} T_1$. Since $z\in D_1
\subset\Omega(G_1)$ and every point of $B^\circ_2\cap\mathrm{Fr}
T_1$ is either a $(G_1-J)$-translate of a point of $S$ or a limit
point of $G_1$, we deduce that $z$ is a $(G_1-J)$-translate of a
point of $S$. On the other hand,  since $z$ is contained in
$C_1=\bar \R^n \setminus T_1$, we see that $z$ is not a
$(G_1-J)$-translate of a point of $S$. This is a contradiction.
\end{proof}

Since $z \in D_1\cap S=D_2\cap S$, as was shown in the proof of
Lemma \ref{key lemma}-(10), no $(G_m-J)$-translates of $B_m$ pass
through $z$ nor accumulate at $z$. Therefore, we have $z\in
C^\circ_1$. Since $\{T_n\}$ is decreasing, the
$(G-J)$-translates of $S$ do not accumulate at $z$, for
$(G-J)$-translates of $S$ accumulate at points in $\bar
T_1$, which is disjoint from $C_1^\circ$. This means that $z$ cannot
be a limit point of $G$; hence $z\in \Omega(G)$. Thus we have shown
that $D$ is contained in $\Omega(G)$.

By Lemma \ref{key lemma}-(6)  and Lemma \ref{precisely
invariant}, we see that $(D_1\cap B^\circ_2)\cup(D_2\cap B^\circ_1)$
is precisely invariant under $\{id\}$ in $G$. Setting $A=(D_1\cap
B^\circ_2)\cup(D_2\cap B^\circ_1)$, we have $D=A\cup(D_1\cap S)$ and
$A \subset C_1^\circ$. Then for any $g\in G-\{id\}$, we have
$g(D)\cap D=(g(A)\cap (D_1\cap S))\cup(g(D_1\cap S)\cap
A)\cup(g(D_1\cap S)\cap(D_1\cap S)).$

 If $g\in J-\{id\}$, then $g(D_1\cap
S)\subset S\setminus D_1$ and $g(A)\cup A\subset B^\circ_1\cup B^\circ_2$.
Therefore, $g(D_1\cap S)\cap(D_1\cap S)=\emptyset$, $g(D_1\cap
S)\cap A=\emptyset$ and $g(A)\cap (D_1\cap S)=\emptyset$. It
follows that $g(D)\cap D=\emptyset$ in this case.

If $g\in G_m-J$, then $g(D_1\cap S)=g(D_m\cap S)\subset T_1$ and
Lemma \ref{key lemma}-(4) and (6) imply that $g(A)\subset
B^\circ_{3-m}$. Since $A\cup(D_1\cap S)=D$ is contained in $C_1$ by
Lemma \ref{in C_1}, and $g(D_1 \cap S)$ is contained in $T_1$, we
have $g(D_1\cap S)\cap A=\emptyset$. We also have  $g(D_1\cap
S)\cap(D_1\cap S)=\emptyset$ since $D_1\cap S= D_2 \cap S$ and $D_1,
D_2$ are fundamental sets of $G_1, G_2$ respectively, and $g(A)\cap
(D_1\cap S)=\emptyset$ since $g(A)$ is contained in $B^\circ_{3-m}$
as was seen above. Therefore also in this case, we have $g(D)\cap
D=\emptyset$.

Now, we consider $g=g_n\circ\cdots\circ g_1\in G-(G_1\cup G_2)$,
where $g_1\in G_m-J$.
Then $g(D_1\cap S)=g(D_m\cap S)\subset g(B_m)\subset
T_n\subset T_1$ and $g(A)=g(D_m\cap B^\circ_{3-m})\cup
g(D_{3-m}\cap B^\circ_m)\subset g_n\circ\cdots\circ
g_2(B^\circ_{3-m})\cup g(B^\circ_m) \text{(Lemma \ref{key lemma}-(6))}
\subset T^\circ_{n-1}\cup
T^\circ_n\subset T^\circ_1\subset B^\circ_1\cup B^\circ_2$.
 These facts imply that $g(D_1\cap
S)\cap(D_1\cap S)=\emptyset$ by Lemma \ref{in C_1},
$g(D_1\cap S)\cap A=\emptyset$ by the fact that $A \subset C_1^\circ$, and
$g(A)\cap (D_1\cap S)=\emptyset$.
Thus we have shown that $D$ is  precisely invariant under $\{id\}$ in $G$.
Since we have already shown that $D\subset\Omega(G)$,
this means that $D\subset$$^\circ\Omega(G)$.
\end{proof}

\begin{lem}
\label{invariant discontinuity}
$S\cap \Omega(J)=S\cap \Omega(G),$ and $S\cap \Omega(J)$ is
precisely invariant under $J$ in $G$.
\end{lem}

\begin{proof}
 Let $z$ be a point in $S\cap\Omega(J)$.
Since $S\cap\Omega(G_m)=S\cap\Omega(J)$ for each $m$ by Lemma
\ref{key lemma}-(2), we have $z\in\Omega(G_m)$. As was
shown in the proof of Lemma \ref{key lemma}-(10), no
$(G_m-J)$-translates of $B_m$ pass through $z$ nor accumulate at
$z$. Therefore $z$ is contained in $C^\circ_1$.

Suppose, seeking a contradiction, that $z$ lies in $\Lambda(G)$.
Then there is a sequence $\{g_k\}$ of
distinct elements of $G$ such that $g_k(y)\rightarrow z$ for all $y$ with at
most one exception.
Since $z$ is contained in $\Omega(G_1)\cap\Omega(G_2)$,
we can assume $|g_k|>1$ for all $k$ by taking a subsequence.
We deduce from the fact that $g_k(B)\subset T_1$ for $B=B_1$ or $B_2$
that $z$ must be contained in $\bar T_1$, which is a contradiction.
Thus we have shown
that $S\cap\Omega(J)$ is contained in $S\cap\Omega(G)$.
The opposite inclusion
is trivial.

Now we turn to prove the latter half of our lemma. It is clear that
$J$ keeps $S\cap\Omega(J)$ invariant. Suppose that there are  points
$y$ and $z$ in $S\cap\Omega(G)=S\cap\Omega(J)$ and that there is an
element $g\in G-J$ such that $g(y)=z$. Express $g$ in a normal form
$g=g_n\circ\cdots\circ g_1$. Then $n>1$ since $S$ is a $(J,
G_m)$-block $(m=1,2)$. Clearly $z$ lies on $g(S)\cap S$.
Moreover since $g(S) = g_n(g_{n-1} \circ \dots \circ g_1(S))$ and
$S$ is contained in both $B_1$ and $B_2$, by Lemma
\ref{interactive}, $g(S)$ is contained in either $g_n(B_m)$, where
$g_n$ is assumed to lie in $G_m$. If $z \in g(S)$ is contained in
$g_n(B_m^\circ)$, then it must lie in $B_{3-m}^\circ$, which
contradicts our assumption. Therefore $z$ must lie in $g_n(S)$. We
may assume that $g_n\in G_1-J$ by interchanging the indices if
necessary. Since $B_1$ is a $(J,G_1)$-block, $B_1 \cap \Omega(G_1)$
is precisely invariant under $J$ in $G_1$, which means that
$g_n(\Omega(G_1) \cap B_1)$ is contained in $B_2^\circ$. Because we
have shown that $z$ lies in $S \cap g_n(S)$, this implies that $z\in
\Lambda(G_1)\subset\Lambda(G)$. Since $z=g(y)\in\Omega(G)$, this is
a contradiction. Thus we have shown that
$g(S\cap\Omega(G))\cap(S\cap\Omega(G))=\emptyset$ for any $g\in
G-J$.
\end{proof}

\begin{proof}[Proof of (3)]
Let $g$ be an element of $G$ which  is not conjugate to any element
of either $G_1$ or $G_2$, such that  $|g|$ is minimal among all
conjugates of $g$ in $G$. Clearly, we have $|g|>1.$ Express $g$
 in a normal form $g=g_n\circ\cdots\circ g_1$. If the
length of $g$ is odd, say, $g_n,\ g_1\in G_m-J$, then $g^{-1}_n\circ
g\circ g_n=g_{n-1}\circ\cdots\circ (g_1\circ g_n).$ The
corresponding normal form of $g^{-1}_n\circ g\circ g_n$ has length
less than $n$, which contradicts the minimality of $|g|$. Therefore
the length of $g$ must be even and $g$ must be a $(3-m, m)$-form.
This implies that $g(B_m)\subset g_n\circ g_{n-1}(B_m)\subset B_m$.
Since $(B^\circ_1, B^\circ_2)$ is a proper interactive pair by
assumption, the last inclusion is proper by Lemma \ref{interactive}.
Hence $g$ has the infinite order and has a fixed point in
$g(B_m)\subset B_m$. Similarly, $g^{-1}(B_{3-m})\subset
g^{-1}_1\circ g^{-1}_2(B_{3-m})\subset B_{3-m},$ where the last
inclusion is proper. Therefore $g$ also has a fixed point
in $g^{-1}(B_{3-m})\subset B_{3-m}$, which may coincide with the
above-mentioned  fixed point.

Since $G$ is discrete and $g$ has infinite order, $g$ is not elliptic.
If $g$ is parabolic, then its fixed point is unique, which we denote by $x$.
Hence the two fixed points mentioned above are equal and $x$ lies on  $S\bigcap g(S)$.
By Lemma \ref{invariant discontinuity}, $x$ is a limit point of $J$.
Since $J$ is geometrically finite, $x$ is either a parabolic fixed point of $J$
or a conical limit point for $J$ by Proposition \ref{limit set}.
Since a conical limit point
for $J$ is also that for $G$ and a conical limit point cannot be a parabolic
fixed point, we see that $x$ is a parabolic fixed point of $J$.
\end{proof}

\begin{proof}[Proof of (4)]
Since $B_1$ and $B_2$ are both blocks, for every parabolic fixed
point $z$ of $J$ with the rank of $Stab_J(z)$ being less than $n$, the peak domain centered at $z$ for $J$ has trivial
intersection with $S=\mathrm{Fr} B_1=\mathrm{Fr} B_2$. This shows
the second condition in the definition of blocks holds for $S$.
Lemma \ref{invariant discontinuity} implies that the first condition
in the definition holds for $S$, hence that $S$ is a
$(J,G)$-block. By Lemma \ref{precisely embedded}, $S$ is precisely
embedded in $G$.
\end{proof}

\begin{proof}[Proof of (5)]
 By (4) shown above, we know that $S$ is a $(J, G)$-block.
 Then (5) follows from Theorem \ref{diameter vanishes}.
\end{proof}

\begin{lem}
\label{C_1 invariant}
$C_1\cap B^\circ_m$ is precisely invariant under $G_{3-m}$ in $G$.
\end{lem}

\begin{proof}
It is obvious that $C_1\cap B^\circ_m=\bar \R^n-G_{3-m}(B_{3-m}).$
Since $G_{3-m}(B_{3-m})$ is invariant under $G_{3-m}$, its complement
$C_1\cap B^\circ_m$ is also invariant under $G_{3-m}$.

If $g\in G_m-J$, then $g(C_1\cap B^\circ_m)\subset
g(B^\circ_m)\subset B^\circ_{3-m},$  and we are done.
Now we consider a general $g$ which is expressed in a normal form
$g=g_n\circ\cdots\circ
g_1$  with $|g|>1$.
If $g$ is an $(m, m)$-form, then
$g(C_1\cap B^\circ_m)\subset g(B^\circ_m)\subset B^\circ_{3-m}$ by
Lemma \ref{interactive}.
If $g$ is an $(m, 3-m)$-form, then $g(C_1\cap B^\circ_m)=
g_n\circ\cdots\circ g_1(C_1\bigcap B^\circ_m)=g_n\circ\cdots\circ
g_2(C_1\cap B^\circ_m)$ as was shown in the last paragraph, and this
last term is contained in $B^\circ_{3-m}$ since
$g_n \circ \dots \circ g_2$ is an $(m,m)$-form.
 If $g=g_n\circ\cdots\circ g_1$ is a $(3-m, k)$-form, where either $k=1$
or $k=2$, then, by the discussion above, we see
$g_{n-1}\circ\cdots\circ g_1(C_1\cap B^\circ_m)\subset
B^\circ_{3-m}$; hence $g(C_1\cap B^\circ_m)\subset
g_n(B^\circ_{3-m})\subset T^\circ_1$.
Thus in every case, if $g\notin
G_{3-m}$, then $g(C_1\cap B^\circ_m)\cap(C_1\cap
B^\circ_m)=\emptyset.$
\end{proof}

\begin{lem}
\label{C}
The set $C$ is contained in the union of  $\Omega(G)\setminus$$^\circ\Omega(G)$
and the $G$-translates of  $D\cup \Lambda(G_1)\cup \Lambda(G_2)$.
\end{lem}

\begin{proof}
Every point $x\in C$  is contained either in $C_1$ or in
$C_n\setminus C_{n-1}$ for some index $n$  $(n>1)$ since $\{C_n\}$
is increasing. If $x\in C_n\setminus C_{n-1}$, then $x\in T_{n-1}\setminus T_n$. Hence
there are a point $y\in B_k$ and  an element expressed in an
$(m,k)$-form $g=g_{n-1}\circ\cdots\circ g_1\in G$ such that
$x=g(y)$. If $y$ lies in $T_1$, then either $y\in (G_k-J)(B_k)\cap
B_k$ or $y\in (G_{3-k}-J)(B_{3-k})$. In the former case, $y$ is
contained in  $\Lambda(G_k)\cap S=\Lambda(J)\cap S$ by Lemma
\ref{key lemma}-(5). In the latter case, we have $x\in T_n$, which
is a contradiction. Therefore, every point $x \in C$ is either
contained in $G(\Lambda(J))$ or $G(C_1)$. In the former case, we are
done. Therefore, we have only to consider the latter case. Moreover,
since the sets in our statement are $G$-invariant, we can assume
that $x$ lies in $C_1$.

It suffices to prove our lemma under the assumption that $x\in
C_1\cap B_2$; the proof for the case $x\in C_1\cap B_1$ is  the
same. If $x$ lies in $C_1 \cap B_2$, then either $x\in \Lambda(G_1)$
or $x\in$$^\circ\Omega(G_1)$ or $x\in
\Omega(G_1)\setminus$$^\circ\Omega(G_1)$. We only need to
discuss the latter two cases.
\smallskip

\noindent{\bf Case $1$:}\quad $x\in$$^\circ\Omega(G_1)$.

In this case, there are an element $g\in G_1$ and a point $z\in D_1$
with $g(z)=x$. We claim that $z\notin B^\circ_1$. Suppose, on the
contrary, that $z$ is contained in $B^\circ_1$. If $g$ lies in $G_1
- J$, then $g(z)$ is contained in $T_1$ by the definition
of $T_1$. Since we assumed that $x$ lies in $C_1$, this is not
possible. Therefore, we have $g \in J$. On the other hand,
$J(B^\circ_1)=B^\circ_1$, which contradicts the assumption that $x$
lies in $B_2$. This shows that $z \in D_1\cap B_2\subset D$,
and we are done in this case.
\smallskip

\noindent{\bf Case $2$:}\quad $x\in \Omega(G_1)\setminus$$^\circ\Omega(G_1)$.

Since $S\cap\Omega(J)=S\cap\Omega(G_1)=S\cap\Omega(G_2)=S\cap
\Omega(G)$ by Lemma \ref{invariant discontinuity}, if $x\in S$, then
$x$ lies in $\Omega(G)$. Furthermore, since $^\circ \Omega(G)$ is
contained in $^\circ \Omega(G_1)$, this implies that
$x\in\Omega(G)\setminus$$^\circ\Omega(G)$, and we are done in this
case. If $x\notin S$, then $x\in C_1\cap B^\circ_2$. Since
$x\in\Omega(G_1)$, no $(G_1-J)$-translates of $B_1$ accumulate at
$x$ as was shown in the proof of Lemma \ref{key lemma}-(10).
Therefore, we have $x\in C^\circ_1$. Then, by Proposition
\ref{domain of discontinuity}, there is a neighbourhood
$U$ of $x$ contained in $C_1\cap B^\circ_2$ such that $U$ is
precisely invariant under $\mathrm{Stab}_{G_1}(x)$ in $G_1$ and
$\mathrm{Stab}_{G_1}(x)$ is a non-trivial finite subgroup. Now Lemma
\ref{C_1 invariant} implies that
$\mathrm{Stab}_{G_1}(x)=\mathrm{Stab}_G(x)$.
Hence $U$ is precisely
invariant under $\mathrm{Stab}_{G}(x)$ in $G$. This shows that $x$
is contained in $\Omega(G)\setminus$$^\circ\Omega(G)$, and we have
completed the proof.
\end{proof}
\begin{lem}
\label{T}
$T\subset\Lambda(G)$.
Furthermore, every point of $T$ is either a $G$-translate of a point in
$\Lambda(J)$ or the limit of nested translates of $S$.
\end{lem}
\begin{proof}
Consider a point $z\in T$. We assume that $z\in (G_1-J)(B_1)$, for
the case when $z \in (G_2-J)(B_2)$ can be dealt with in the same
way. Then there is an element $h_1=g_1\in G_1-J$ such that $z\in
g_1(B_1)$. Since $T_1 \supset T_2$, we have $z\in T_2$, and there is
an element $g_2\in G_2-J$ such that $z\in g_1\circ
g_2(B_2)=h_2(B_2)\subset h_1(B_1)$. Similarly, since $z\in T_3$,
there is an element $g_3\in G_1-J$ such that $z\in g_1\circ g_2\circ
g_3(B_1)=h_3(B_1)\subset h_2(B_2)\subset h_1(B_1)$; etc. Since the
element $h_k$ has length increasing as $k \rightarrow \infty$ and
$(B^\circ_1, B^\circ_2)$ is a proper interactive pair,  the sets
$h_k(S)$ can be assumed to be all distinct by taking a subsequence
if necessary. Thus we have shown that if $z\in T$, then there is a
sequence $\{h_k\}$ of elements of $G$, with
$|h_k|\rightarrow\infty$, and $z\in\cdots\subset
h_k(\check{B}_k)\subset\cdots\subset h_2(\check{B}_2)\subset
h_1(\check{B}_1)$, where $\check{B}_j$ is either $B_1$ or $B_2$.
Passing to a subsequence if necessary, we may assume that
$\check{B}_j=B_1$.

There are two possibilities for this sequence: either $z$ lies in
the interiors of infinitely many  $h_k(B_1)$, or from some
$k$ on, $z$ lies on the boundary of every $h_k(B_1)$. In
either case, since the $h_k(S)$ are distinct, we have
$\mathrm{diam}(h_k(S)) \rightarrow 0$. Since the ball
$h_k(B_1)$ bounded by $h_k(S)$ decreases as $k \rightarrow
\infty$, this is possible only when $\mathrm{diam}(h_k(B_1))
\rightarrow 0$. Since $z$ is a limit of $\{h_k(x_k)\}$ with $x_k \in
B_1$ in either case above, it follows that for every $x
\in B_1$, we have $h_k(x) \rightarrow z$. This means that
$z$ lies in $\Lambda(G)$. Moreover, in the former case, we have
shown that  $\{h_k(S)\}$ nests around $z$. In the latter case, since
$z\in h_{k_0}(S)\cap h_{k_0+1}(S)\cap\cdots$, we have
$w=h^{-1}_{k_0}(z)\in S\cap h^{-1}_{k_0}h_{k_0+1}(S)\cap\cdots$.
Since such $w$ is contained in $\Lambda(G)$, by Lemma \ref{invariant
discontinuity}, it also lies in $\Lambda(J)$. This means that $z$ is
contained in the $G$-translate of $\Lambda(J)$.
%
%
\end{proof}

\begin{lem}
\label{nesting}
If $z\in C\cap\Lambda(G)$, then there is no sequence of distinct
translates of $S$ nesting about $z$.
\end{lem}

\begin{proof}
 Lemma \ref{C} implies that  $z$ is a $G$-translate of a point
 in either $D$ or $\Lambda(G_1) \cup \Lambda(G_2)$.
 Since $D$ is contained in $\Omega(G)$ by Lemma \ref{D}, the only
 possibility is $z \in G(\Lambda(G_1) \cup \Lambda(G_2))$.

We first consider the special case when $z$ lies in $G(\Lambda(J))$.
Under this assumption, suppose, seeking a contradiction, that there
is a sequence $\{h_k(S)\}$ of distinct $G$-translates of $S$ nesting
about $z= g(y)$ for an element $g\in G$ and a point $y\in
\Lambda(J)\subset S$. Then we have $z\in h_k(B^\circ)$ by taking a
subsequence for $B$ which is either $B_1$ or $B_2$. We can assume
that $B$ is $B_1$ after taking a subsequence, for we can deal with
the other case in the same way. It follows that $y\in g^{-1}\circ
h_k(B_1^\circ)$. Now since $\{h_k(S)\}$ nests around $z$, we have
$\mathrm{diam}(h_k(B_1)) \rightarrow 0$. This is possible only when
after taking a subsequence all $h_k$ are $(m_k,1)$-forms with
$m_k=1,2$. (If $h_k$ were $(m_k ,2)$-form, then $h_k(B_1)$ would contain
$S$; hence its diameter would not go to $0$.) Therefore $g^{-1}h_k$
is also expressed as an $(m',1)$-form for large $k$ and
$g^{-1}h_k(B_1^\circ)$ is contained in $B_{3- m'}^\circ$.
In particular, we have $y\notin S$. This contradiction shows that if
$z\in G(\Lambda(J))$, then there is no sequence of distinct
translates of $S$ nesting about $z$.

Now we turn to the general case when $z \in G(\Lambda(G_1) \cup \Lambda(G_2))$.
It suffices to consider the case $z\in G(\Lambda(G_1))$ since the
proof for the case $z\in G(\Lambda(G_2))$ is entirely the same.
Then there are an element $g\in G$ and  a point $y\in\Lambda(G_1)$ with $g(y)=z$.
Since $B^\circ_1\subset\Omega(G_1)$,
we have $\Lambda(G_1)\subset \bar \R^n \setminus G_1(B^\circ_1)$.
Therefore, $y$ is not contained in $G_1(B^\circ_1)$; hence unless $y$
lies in $G_1(S)$, it must lie in $C_1 \cap B^\circ_2$.
 If $y\in G_1(S)$, then $y\in
G_1(S\cap\Lambda(G_1))=G_1(S\cap\Lambda(J))$.
The discussion in the previous paragraph implies that this case cannot occur.

Now we assume that $y\in C_1\cap B^\circ_2$.
If there is a sequence $\{h_k(S)\}$ of
distinct $G$-translates of $S$ nesting about $z= g(y)$, then $z\in
h_k(B^\circ)$ for every $k$ where $B$ is $B_1$ or $B_2$, and
hence $y\in g^{-1}\circ h_k(B^\circ)$.
We may assume that $B=B_1$ by changing the index and taking a
subsequence and $h_k$ is an $(m,1)$-form.
Then $g^{-1}\circ h_k$ is also an $(m',1)$-form for sufficiently large $k$.
Since $\{T_n\}$ is a decreasing sequence,
$y\in T^\circ_1$, which is a contradiction.
Thus we have completed the proof.
\end{proof}

\begin{proof}[Proof of (6)]
 If $x$ lies in $\Lambda(G)\setminus G(\Lambda(G_1)\cup\Lambda(G_2))$,
 then $x\in T$ by Lemma \ref{C}.
 Since every point of $T$ is either a translate of a
point of $\Lambda(J)$ or is the limit of a nested sequence of
translates of $S$ by Lemma \ref{T}, we have proved the ``if" part.

Now we turn to the ``only if" part.
Suppose that $x$ lies in $\Lambda(G_m)$ for $m=1$ or $2$.
Since  $B^\circ_m\subset\Omega(G_m)$ by Lemma \ref{key lemma}-(3),
we have $x\in \bar \R^n \setminus G_m(B^\circ_m)$.
If $x\in G_m(S)$, then as was shown in  the proof of
Lemma \ref{nesting}, there is no
distinct $G$-translates of $S$ nesting about $x$.
Therefore $x$ is contained in $\bar \R^n \setminus G_m(B_m)=C_1\cap
B^\circ_{3-m}$, which implies that $x \in C \cap \Lambda(G)$.
By Lemma \ref{nesting}, there is no distinct translates of $S$
nesting about $x$.
\end{proof}

\begin{proof}[Proof of (7)]
By Lemma \ref{C}, every point of
$C\cap$$^\circ\Omega(G)$ is a translate of a point of $D$.
Also by Lemma \ref{T}, $T$ is contained in $\Lambda(G)$.
This shows that every point of $^\circ \Omega(G)$ is contained in
a $G$-translate of $D$.
Furthermore, since
$D\subset$$^\circ\Omega(G)$ and $D$ is precisely invariant under
the identity in $G$ by Lemma \ref{D}, it follows that $D$ is
a fundamental set for $G$.

Now assume that both $D_1$ and $D_2$ are constrained.

\begin{cl}
 $\Omega(G)\subset G(\bar D)$.
 \end{cl}
\begin{proof}
Since we have already shown that $D$ is a fundamental set
for $G$, we have only to prove that if
$x\in\Omega(G)\setminus$$^\circ\Omega(G)$, then there is an element
$g\in G$ with $g(x)\in\bar D$. Now let $x$ be a point in
$\Omega(G)\setminus$$^\circ\Omega(G)$. By Lemma \ref{T}, $x$ is not
contained in $T$. As was shown in the proof of Lemma \ref{C}, we
have $x\in G(C_1)\cap(\Omega(G)\setminus$$^\circ\Omega(G))$. This
means that there are an element $g\in G$ and a point $y\in
C_1\cap(\Omega(G)\setminus$ $^\circ\Omega(G))$ such that $x=g(y)$.
We may assume that $y\in B_2$, for the proof in the case $y\in B_1$
is entirely the same.

Suppose first that  $y\in S \cap
C_1\cap(\Omega(G)\setminus$$^\circ\Omega(G))$. Then since
$S\cap\Omega(J)=S\cap\Omega(G_1)=S\cap\Omega(G)$ by Lemma
\ref{invariant discontinuity} and $D_1$ is a constrained
fundamental set for $G_1$, there are an element $h\in G_1$ and a
point $z\in\bar D_1$ such that $y=h(z)$. Since
$G_1(B^\circ_1)\subset B^\circ_1\cup T^\circ_1$, we see that $z$
must be contained in $B_2$, hence $z\in \bar D_1\cap B_2\subset \bar
D$. Thus we have completed the proof in this case.

Next we assume that $y\notin S$, which means that
$y\in C_1\cap B^\circ_2\cap(\Omega(G)\setminus $$^\circ\Omega(G))$.
Since $y\in\Omega(G)\subset\Omega(G_1)$ and $D_1$
is a fundamental set for $G_1$, we see that $y$ is
$G_1$-equivalent to a point $w\in\bar D_1$.
By Lemma \ref{C_1 invariant}, we have $w\in\bar D_1\cap C_1\cap B^\circ_2$.
Since $\bar D_1\cap
B^\circ_2\subset\bar{D}$, this implies $w\in \bar D$, and  our claim has been proved.
\end{proof}

We now return to the proof of (7).
We have
\begin{equation}
\label{first}
G_m(\bar D_m)=G_m((\bar D_m\cap
B^\circ_m)\cup (\bar D_m\cap B_{3-m})),
\end{equation}
\begin{equation}
\label{second}
G_m(\bar
D_m\cap B^\circ_m)\subset B^\circ_m\cup (T^\circ_1\cap
B^\circ_{3-m})
\end{equation}
 by the definition of $T_1$, and
\begin{equation}
\label{third}
 \bar D_m\cap
B^\circ_{3-m}\subset\overline{D_m\cap B_{3-m}}\subset \bar C_1\cap
B_{3-m}
\end{equation}
by Lemma \ref{in C_1}.

Since $\bar C_1\cap B_{3-m}=\bar \R^n \setminus  G_m(B^\circ_m)$, we
see that $\bar C_1\cap B_{3-m}$ is $G_m$-invariant.
Therefore from (\ref{third}), we obtain
\begin{equation}
\label{fourth}
G_m(\bar D_m \cap B^\circ_{3-m}) \subset \bar C_1 \cap B_{3-m}.
\end{equation}
Since $\mathrm{Fr} D\cap S$ consists of
only finitely many connected components the sum of whose
$(n-1)$-dimensional measures on $S$ vanishes by assumption, it follows
from (\ref{first}), (\ref{second}), and (\ref{fourth})  that the sides of
$D_m$ in $B_{3-m}$ are paired with those in $B_{3-m}$ by elements of $G_m$
for each $m$. Since the sides of $D$ in $B_1$ are equal to those of $D_2$ in
$B_1$ and the sides of $D$ in $B_2$ those of $D_1$ in $B_2$, we
see the sides are paired to each other.
These sides can accumulate only at limit points because of the same
property for $D_1$ and $D_2$.
The only thing left to show is that the tessellation of $\Omega(G)$
by translates of $\bar D$ is locally
finite.

Take any $z\in \bar D\cap\Omega(G)$. We see from Lemma \ref{in C_1}
that either $z\in C^\circ_1$ or $z\in\mathrm{Fr}
C_1=\mathrm{Fr} T_1$. We may assume that $z\in
\overline{D_1\cap B_2}\subset \bar D_1\cap B_2$, for the proof in
the case $z\in \bar D_2\cap B_1$ is entirely the same.
\smallskip

\noindent{\bf Case 1}\quad $z\in C^\circ_1$.

Since $z$ is contained in $\Omega(G_m)$ for each $m$ and $D_m$ is a
constrained fundamental set for $G_m$, there is a
neighborhood $U$ of $z$ with $U\subset C_1^\circ$ such that for each
$m$ there is a finite set $\{g_{m1}(D_m),\cdots, g_{m{k_m}}(D_m)\}$
with $U\subset\cup_i g_{mi}(\bar D_m)$ for $g_{mi}\in G_m$. We
consider $U\cap B_{3-m}$. Since $G_m(\bar D_m\cap B^\circ_m)\subset
B^\circ_m\cup T^\circ_1$ and $U\subset C_1$, we have $U\cap
B_{3-m}\subset\cup_i g_{mi}(\bar D_m\cap B_{3-m}).$ Therefore
$U\subset\cup^2_{m=1}(\cup g_{mi}(\bar D_m\cap
B_{3-m}))\subset\cup^2_{m=1}(\cup_i g_{mi}(\bar D)),$ and we have
obtained the local finiteness of $D$  at such a point.
%

\smallskip
\noindent{\bf Case 2}\quad $z\in\mathrm{Fr} C_1=\mathrm{Fr} T_1$.

We claim that $z\notin S$ in this case. Suppose, on the contrary,
that $z$ is contained in $S$. Since $z\in
\Omega(G)\subset\Omega(G_m)$, as was shown in the proof of Lemma
\ref{key lemma}-(10), no $(G_m-J)$-translates of $B_m$
pass through $z$ and no $G_m$-translates of $B_m$ accumulate at $z$.
Therefore, we have $z\in C^\circ_1$, which contradicts our
assumption for Case 2.

Hence, we can assume that  $z$ lies in $B^\circ_2$. Since a point of
$\mathrm{Fr} T_1$ in $B^\circ_2$ is either a point of $(G_1-J)(S)$,
or a point of $\Lambda(G_1)$ and $z\in\Omega(G)\subset\Omega(G_1)$,
we see that $z$ must lie in $B^\circ_2\cap(G_1-J)(S)$. Then there
are a point $s\in S$ and an element $g\in G_1-J$ with $g(s)=z$. By
Lemma \ref{invariant discontinuity}, $s$ lies in $S\cap
\Omega(G)=S\cap \Omega(J)=S\cap\Omega(G_1)=S\cap\Omega(G_2)$.
Therefore no $(G_m-J)$-translates of $B_m$ pass through $s$ and no
$G_m$-translates of $B_m$ accumulate at $s$ as was shown in the
proof of Lemma \ref{key lemma}-(10). This implies that $s$ is
contained in $C^\circ_1\cap S$. By applying the proof of Case
$1$ to $s$, we see that  there is a neighbourhood $U$ of
$s$ covered by finitely many $G$-translates of $\bar D$. It follows
that  $g(U)$ is a neighbourhood of $z$ covered by finitely many
$G$-translates of $\bar D$. This shows that  $D$ is locally finite
at a point as in  Case 2.

Thus we have shown the proof of the local finiteness of $D$, hence completed the proof.
\end{proof}

\begin{proof}[Proof of (8)]
We shall prove this  by showing the following three claims.

\begin{cl}
\label{R_m and Omega}
 For each $m$, we have
$R_m\cap\Omega(G_m)\subset\Omega(G)$.
\end{cl}
\begin{proof}
Take a point $z\in R_m\cap\Omega(G_m)$. Since $R_m=\bar \R^n
\setminus G_m(B^\circ_m)$, we have either $z\in G_m(S)$ or $z\in
C_1\cap B^\circ_{3-m}$. If $z\in G_m(S)$, then $z\in\Omega(G)$ since
$S\cap\Omega(G)=S\cap\Omega(J)=S\cap\Omega(G_m)$ by Lemma
\ref{invariant discontinuity}. If $z\in C_1\cap
B^\circ_{3-m}$, since $z\in\Omega(G_m)$, no $G_m$-translates of
$B_m$ passe through or accumulate at $z$ as was shown in the proof
of Lemma \ref{key lemma}-(10). It follows that $z\in C^\circ_1$. By
Proposition \ref{domain of discontinuity}, there is a neighbourhood
$U$ of $z$ lying in $C_1^\circ \cap B^\circ_{3-m}$ which is
precisely invariant under $\mathrm{Stab}_{G_m}(z)$ in $G_m$ such
that $\mathrm{Stab}_{G_m}(z)$ is finite. By Lemma \ref{C_1
invariant}, we see that $\mathrm{Stab}_{G_m}(z)=\mathrm{Stab}_G(z)$
and that $U$ is precisely invariant under $\mathrm{Stab}_{G}(z)$ in
$G$. By Proposition \ref{domain of discontinuity}, this implies that
$z\in\Omega(G)$.
\end{proof}

\begin{cl}
Every point of $\Omega(G)$ is $G$-equivalent to a point of either
$R_1\cap\Omega(G_1)$ or $R_2\cap\Omega(G_2)$.
\end{cl}

\begin{proof}
Let $z$ be a point in $\Omega(G)$.
By Lemma \ref{T}, we see that $z\notin T$.
As was shown in the first half of the proof of Lemma \ref{C}, we have $z\in G(C_1)$.
We have only to consider the case when $z\in C_1$ by translating $z$ by elements of $G$.
Since $C_1\cap B_m\subset R_{3-m}$ by the definitions of $R_{3-m}$ and $C_1$ and
$\Omega(G)\subset
\Omega(G_1)\cap \Omega(G_2)$, we see that $z\in (R_1\cap \Omega(G_1))\cup
(R_2\cap \Omega(G_2))$.
\end{proof}

\begin{cl} For each $m=1,2$, the set $R_m\cap\Omega(G_m)$ is
precisely invariant under $G_m$ in $G$.
\end{cl}
\begin{proof}
It is obvious that $R_m$ is $G_m$-invariant, hence so is $R_m\cap\Omega(G_m)$.
We shall show that $R_m \cap \Omega(G_m)$ is moved to a set disjoint from it
by other elements of $G$.

 For any $g\in G_{3-m}-J$, we have
 $g(R_m\cap\Omega(G_m))\subset g(B_{3-m}\cap\Omega(G_m))\subset B_m$.
By Lemma \ref{key lemma}-(5), $g(B_{3-m})\cap S\subset
\Lambda(G_{3-m})\cap S$, which is equal to $S\cap\Lambda(G_m)$ by
Lemma \ref{key lemma}-(2). This implies that no point of
$\Omega(G_m)\cap B_{3-m}$ is mapped into $S$ by $g$, hence
$g(B_{3-m}\cap\Omega(G_m))\subset B^\circ_m$. Since $R_m$ is
contained in $B_{3-m}$, it follows that $g(R_m \cap
\Omega(G_m)) \cap R_m \cap \Omega(G_m) = \emptyset$.

 Now let $g=g_n\circ\cdots\circ g_1$ be a normal form with $|g|>1$.
 If $g$ is a $(3-m, 3-m)$-form, then since $g_1(R_m\cap\Omega(G_m))\subset
B^\circ_m$, we have $g(R_m\cap\Omega(G_m))\subset
g_n\circ\cdots\circ g_2(B^\circ_m)\subset B^\circ_m$.
If $g$ is a
$(3-m, m)$-form, then since $g_1$ preserves $R_m\cap\Omega(G_m)$, we
have $g(R_m\cap\Omega(G_m))=g_n\circ\cdots\circ
g_2(R_m\cap\Omega(G_m))$, which is contained in $B^\circ_m$ by
the argument above for $(3-m,3-m)$-forms.
Finally if $g$ is an $(m, k)$-form, then $g_{n-1}\circ\cdots\circ g_1$
is a $(3-m, k)$-form with $k=3-m$ or $k=m$.
Then, as was discussed above, we have
$g_{n-1}\circ\cdots\circ g_1(R_m\cap\Omega(G_m))\subset B^\circ_m$,
and $g(R_m\cap\Omega(G_m))\subset g_n(B^\circ_m)$, which is
contained in the complement of $R_m$ by definition.
Thus we have shown that
$g(R_m \cap \Omega(G_m))\cap R_m \cap \Omega(G_m)= \emptyset$ for any $g \in G-G_m$.
\end{proof}

By these three claims, we have shown that $\Omega(G)/G=(R_1 \cap
\Omega(G_1))/G_1 \cup (R_2 \cap \Omega(G_2))/G_2$. Now we consider
the intersection of the two terms in the right hand side. We should
first note that $(R_1 \cap \Omega(G_1)) \cap (R_2 \cap \Omega(G_2))$
is contained in $B_2 \cap B_1 =S$ since $R_1$ is contained in $B_2$,
and $R_2$ is in $B_1$. Since $\Omega(G_m) \cap S=\Omega(J)
\cap S\subset R_m\cap \Omega(G_m)$, the intersection is
equal to $\Omega(J) \cap S$. Furthermore since $S$ is a
$(J,G_m)$-block, $\Omega(J)\cap S$ projects to $(\Omega(J) \cap
S)/J$ in $(R_m \cap \Omega(G_m))/G_m$. Therefore $(R_1 \cap
\Omega(G_1))/G_1$ and $(R_2 \cap \Omega(G_2))/G_2$ are pasted along
$(S \cap \Omega(J))/J$.
\end{proof}

In the following,  we assume further that each $B_m$ is precisely
invariant under $J$ in $G_m$.

\begin{proof}[Proof of (9)]
 Let $x$ be a parabolic fixed point of $J$.
 Such a point $x$ is contained in $S$ by Lemma \ref{key lemma}-(2).
 Since each $B_m$ is precisely invariant under $J$ in $G_m$ by our assumption,
we have $\mathrm{Stab}_J(x)=\mathrm{Stab}_{G_m}(x)$, which is also
equal to $\mathrm{Stab}_G(x)$ by Lemma \ref{precisely
embedded}. Let $H$ denote $\mathrm{Stab}_J(x)$.

\noindent{\bf The ``if" part.}\quad Suppose that $B_m$ is a strong
$(J,G_m)$-block for each $m=1,2$. There is nothing to prove if the
rank of $H$ is $n$, for the rank of $\mathrm{Stab}_G(x)$ is also $n$
then. Now assume that the rank of $H$ is $k<n$. By conjugation, we
may assume that $x=\infty$. By Theorem \ref{Bieberbach}, we can
assume that each $g\in H$ is expressed as $g(x)=Ax+\mathbf{a}$ for
$\mathbf{a}\in \R^k$ and an orthogonal matrix $A$
preserving the subspaces $\R^k$ and $\R^{n-k}$.

Since both $B_1$ and $B_2$ are assumed to be strong and
$\mathrm{Stab}_{G_1}(\infty)=\mathrm{Stab}_{G_2}(\infty)$, there is a common peak domain $U$ at $\infty$
for $G_1$ and $G_2$. Since
$U\cap(\Lambda(G_1)\cup\Lambda(G_2))=\emptyset$, by choosing $U$
small enough, we may assume that $\bar U\setminus
\{\infty\}\subset\Omega(G_1)\cap\Omega(G_2)$, where $\bar{} $ means
the closure on $\bar \R^n$. We can assume that $U$ has a form
$U=\{x\in R^n: \sum^n_{i=k+1}x^2_i>t^2\},$ where $t$ is a
sufficiently large positive number.

\begin{cl}
\label{small peak}
We can choose $U$ small enough to satisfy $U\subset
C_1$.
\end{cl}

\begin{proof}[Proof of Claim]
We divide our discussions into two cases.
\smallskip

\noindent{\bf Case 1}\quad The case when $k=n-1$.

In this case,  $U$ is the union of two components $U_1$ and $U_2$, and
we may assume that $U_m\subset B^\circ_m$ by our assumption that $B_m$
is a strong block.
We have only to prove that we can choose $U_1$ small enough in such a
way that every $G_2$-translate of $B_2$ is disjoint from $U_1$.
We may assume that $U_1=\{x\in \R^n: x_n>t\}$.
Suppose, seeking a contradiction, that such a $U_1$ does not exist.
Then, there is a sequence $\{g_k(B_2)\}$ of distinct
$G_2$-translates of $B_2$  intersecting
$\{x\in \R^n: x_n>s\}$ for any large $s$.
This means that the projections of $g_k(B_2)$ to the $n$-th coordinate
$\R$ accumulate at $\infty$.
We may assume that $g_k\in G_2-J$ since $J$ fixes $B_2$.

Now Lemma \ref{key lemma}-(7) implies that $\mathrm{diam}(g_k(B_2))\rightarrow 0$
with respect to the ordinary spherical metric.
It follows that  $g_k(y)\to \infty$ for all $y\in B_2$ since $\{g_k(B_2)\}$
accumulates at $\infty$.
By Lemma \ref{converge to limit points}, by taking a  subsequence of $\{g_k\}$, we may
assume that $g_k(y)\rightarrow \infty$ for all $y$ with at most one exception,
which must be a limit point.

 Since $\bar U_2\setminus\{\infty\}$ is contained in
$\Omega(G_2)$, for all $y\in \bar U_2\setminus \{\infty\}$, we have
$g_k(y)\rightarrow \infty$. Since $g_k(U_2)\cap U=\emptyset$, it
follows that the projections of $g_k(\bar U_2)$ to the $n$-th
coordinate are bounded. Hence the projections of $g_k(\bar
U_2\backslash \infty)$ to the first $n-1$ coordinates $\R^{n-1}$
accumulate at $\infty$. By Theorem \ref{Bieberbach}, for each $g_k$,
we can choose an element $j_k\in H$ such that $\{j_kg_k(y_0)\}$ lies
in a bounded set for a fixed $y_0\in U_2$. For each $k$, we have
$\infty\notin g_k(B_2)$ since $B_2$ was assumed to be precisely
invariant under $J$ in $G_2$ and $\infty$ lies on $S$. Therefore, we
have $\infty\notin j_kg_k(B_2)$. Since
$|(j_kg_k(y))_n|=|(g_k(y))_n|$ and the projections of the $g_k(B_2)$
to the $n$-th  coordinate $\R$ accumulate at $\infty$, we see that
$\{j_kg_k(B_2)\}$ also accumulates at $\infty$. By Lemma \ref{key
lemma}-(7), this implies that $j_kg_k(y)\rightarrow \infty$ for all
$y\in B_2$. This is a contradiction since  $\{j_kg_k(y_0)\}$ stays
in a compact set. This proves our claim for the case when $k=n-1$.
\smallskip

\noindent{\bf Case 2}\quad The case when $k<n-1$.

Since $U$ is connected and is disjoint from $S$, we see that $U$
lies in either $B_1^\circ$ or $B_2^\circ$. We may assume that
$U\subset B^\circ_1$. Then, by the same argument as in the proof of
Case 1, we see that the projections of $G_2$-translates of
$B_2$ in the last $n-k$ coordinates cannot accumulate at $\infty$.
Therefore, we have $U\subset C_1\cap B^\circ_1$.
\smallskip

The claim has thus been proved.
\end{proof}

Now we return to the proof of the``if" part of (9). Take a small
common peak domain $U$ for both $G_1$ and $G_2$ as in Claim
\ref{small peak}. By assumption, $U$ is precisely invariant under
$H$ in both $G_1$ and $G_2$. We need to show it is precisely
invariant under $\mathrm{Stab}_G(x)$ in $G$.

For any $g\in G-(G_1\cup G_2)$, we have $g(U)=g(U_1)\cup
g(U_2)\subset g(C_1\cap B^\circ_1)\cup g(C_1\cap B^\circ_2)$, where
$U_1, U_2$ are the components of $U$ if $k=n-1$, and we regard one
of them as the emptyset when $k < n-1$. Suppose that $g$ is
expressed as a $(1,1)$-form $g_n \circ\dots\circ g_1$. As
was shown in Lemma \ref{interactive}, $g_n \circ\dots\circ
g_1(C_1 \cap B^\circ_1) \subset B^\circ_2$. Furthermore, we have
$g_n \circ\dots\circ g_1(C_1 \cap B^\circ_1) \subset g_n
\circ\dots\circ g_1(B^\circ_1) \subset T_n^\circ \subset
T_1^\circ$. On the other hand, $g_n \circ\dots\circ g_1(C_1
\cap B^\circ_2) \subset g_n \circ\dots\circ g_2(C_1 \cap
B^\circ_2)$ by Lemma \ref{C_1 invariant}. Then applying the same
argument for $C_1 \cap B^\circ_1$, we see that $g_n
\circ\dots\circ g_2(C_1 \cap B^\circ_2) \subset T_1^\circ$.
Thus we have shown that $g(C_1 \cap B^\circ_1) \cup g(C_1\cap
B^\circ_2) \subset B_2^\circ \cap T^\circ_1$ for $g$ expressed as a
$(1,1)$-form. A similar argument works also for $(1,2)$-form. Also,
we can see by the same argument that if $g$ is expressed as a
$2$-form, then $g(U)=g(U_1)\cup g(U_2)\subset g(C_1\cap
B^\circ_1)\cup g(C_1\cap B^\circ_2)\subset B^\circ_1\cap T^\circ_1$.

Since $U$, which is disjoint from $S$ from the beginning, is taken
to be lie inside $C_1$,  it follows that $U$ is precisely invariant
under $H$ in $G$ in the case when $k\leq  n-1$.

This completes the proof of the ``if" part.

\noindent{\bf The ``only if" part.}\quad Let $x$ be a parabolic
fixed point of $J$ such that $\mathrm{Stab}_J(x)$ has rank
less than $n$. This point $x$ must lie on $S$ since $\Lambda(J)
\subset S$. Since we are assuming that $S$ is a strong
$(J,G)$-block, there is a peak domain $U$ for $G$, which is also a
peak domain for both $G_1$ and $G_2$. Since we already know that
$B_m$ is a $(J,G_m)$-block, this shows that $B_m$ is a
strong $(J, G_m)$-block.
\end{proof}

\begin{proof}[Proof of (10)]
Since we are assuming both $B_1$ and $B_2$ are strong blocks, by (9), $S$ is a strong $(J,G)$-block.
Let $x$ be a limit point of $G$ which is not a translate of a limit point of either $G_1$ or $G_2$.
By Lemma \ref{C}, we see that $x$ is contained in $T$.
Furthermore, by Lemma \ref{T},  there is a sequence $\{h_k\}$ of distinct elements
of $G$  such that $x\in\cdots h_k(B)\subset\cdots\subset
h_1(B)$ for $B$ which is either $B_1$ or $B_2$.
We can assume that $B=B_1$ and $h_1=id$ by interchanging the indices and  replacing
$g(B_2)$ with $B_1$ for $g \in G_2$ if necessary.
Then $S$ separates $h^{-1}_k(S)$ from $h^{-1}_k(x)$.

Since $J$ is geometrically finite, by Proposition \ref{essentially
finite}, there are a Dirichlet domain $P$ and standard parabolic
regions $B_{p_1}, \dots B_{p_k}$ such that $\bar P \setminus \cup_j
(\mathrm{Int}B_{p_j}\cup \{p_j\})$ is compact. Since $P$
is a Dirichlet domain, the interior of $D=\bar P \cap \bar \R^n$ is
a fundamental domain for $J$. Since $h^{-1}_k(x)\in \Omega(J)$ for
each $k$, there is an element $j_k\in J$ such that $j_k\circ
h^{-1}_k(x)\in D$. We denote $j_k \circ h^{-1}_k$ by $l_k$.

 We claim that $\{l_k(x)\}$ stays away from $S$. Suppose,
 on the contrary, that $l_k(x)\rightarrow w\in S$.
 Then, by Lemma \ref{invariant discontinuity}, $w$ is a parabolic fixed point of $J$,
 where the rank of $\mathrm{Stab}_J(w)$
 is less than $k$ since $D$ intersects $\Lambda(J)$ only at the $p_j$.

This means that all the $l_k(x)$ lie in some $B_{p_j}$ if we
take a subsequence, where $p_j=w$. By the proof of (9), we
can assume that the interior of $B_{p_j} \cap \bar \R^n$, which is
denoted by $U_{p_j}$, is also a peak domain for $G$. Hence
we may assume that $\bar U_{p_j}\setminus
\{p_j\}$ is contained in $\Omega(G)$. On the other hand,
since $x$ lies in $\Lambda(G)$, we have $l_k(x) \in \Lambda(G)$,
which is a contradiction.

Therefore, there is  $\delta>0$ such that $d(l_k(x),z) >\delta$ for all $z \in S$,
where $d$ denotes the ordinary spherical metric on $\bar \R^n$.
Since $S$ separates $h_k^{-1}(x)$ from $h_k^{-1}(S)$, we see that for all $z$ on
$S$ we have $\delta< d(l_k(x), z) \leq d(l_k(x), l_k(z))$.
On the other hand, since $h_k(S)$ nest around $x$, we see that for any point $y$
on $S$, the points $l_k^{-1}(y)$ converge to $x$.
We can now apply Proposition \ref{conical limit} to conclude that $x$ is a conical limit point.
\end{proof}

\begin{proof}[Proof of (11)]
We first assume that $G_1$ and $G_2$ are geometrically finite. Then
every parabolic fixed point of $G_m$ is a parabolic vertex by
Proposition \ref{limit set}. Therefore $B_1$ and $B_2$ are both
strong blocks. By (9), this implies that $S$ is a strong
$(J,G)$-block.

Let $x$ be a point on $\Lambda(G)$.
What we have to show is that $x$ is either a parabolic vertex or a conical
limit point, for this proves that $G$ is geometrically finite by Proposition \ref{limit set}.
Suppose first that $x$ is a parabolic fixed point, where the rank $k$ of
$H=\mathrm{Stab}_G(x)$ is less than $n$.
We shall show that $x$ is a parabolic vertex then.
Since $x$ is a parabolic fixed point, it cannot be a conical limit point.
Hence by (10),  $x$ is a translate of a limit point of either $G_1$ or
$G_2$.

By interchanging the indices and translating $x$ by elements of $G$, we may
assume that $x$ lies in $\Lambda(G_1)$.
Since $G_1$ is assumed to be geometrically finite, $x$ is a
parabolic vertex or a conical limit point for $G_1$ by
Proposition \ref{limit set}.
If $x$ is a conical limit point for $G_1$, then so is it for $G$, which
contradicts the assumption that $x$ is a parabolic fixed point.
Therefore, $x$ is a parabolic vertex for $G_1$.
Suppose first that  $x$ lies on $G_1(S)$.
Then there is an element $\gamma \in G_1$ such that $\gamma^{-1}(x)$ lies on $S$.
Since $x$ is not a conical limit point for $G_1$, neither is $\gamma^{-1}(x)$.
This also implies that $\gamma^{-1}(x)$ is not a conical limit point for $J$ either.
Since $J$ is geometrically finite, again by Proposition \ref{limit set},
we see that $\gamma^{-1}(x)$ is a parabolic vertex for $J$.
Since $S$ is a strong $(J,G)$-block, it follows that $\gamma^{-1}(x)$ is
a parabolic vertex also for $G$, hence so is $x$.
Thus we are done for this case.

Suppose next that  $x$ does not lie on any $G_1$-translate of $S$.
We shall show that $x$ is a parabolic vertex for $G$ even in this
case. Since $G_1(B^\circ_1) \subset\Omega(G_1)$ by Lemma \ref{key
lemma}-(3) and $x$ is a parabolic vertex of $G_1$, we have $x\in
B^\circ_2\cap C_1$. Since $B^\circ_2\cap C_1$ is precisely invariant
under $G_1$ in $G$ by Lemma \ref{C_1 invariant},
$H=\mathrm{Stab}_G(x)$ must be contained in $G_1$. This implies that
$H=\mathrm{Stab}_{G_1}(x)$. Since $x$ is a parabolic vertex for
$G_1$, there is a peak domain $U$ at $x$ for $G_1$. Since
$U\cap\Lambda(G_1)=\emptyset$ and $x\in B^\circ_2\cap C_1$, by
choosing $U$ to be sufficiently small, we can assume that $\bar
U\setminus\{x\}\subset\Omega(G_1)$ and $\bar U\subset B^\circ_2$. By
conjugating $G$ by an element of $M(\bar \R^n)$, we may assume that
$x=\infty$ and $U$ is in the form $U=\{x\in R^n:
\sum^n_{i=k+1}x^2_i>t\},$ for some $t>0$. By Theorem
\ref{Bieberbach}, for any $g \in \mathrm{Stab}_G(\infty)$,
we have an expression $g(x)=Ax+\mathbf{a}$, for $\mathbf{a}\in \R^k$
and an orthogonal matrix $A$ preserving the subspaces $\R^k$ and
$\R^{n-k}$. Now we shall show the following.

\begin{cl}
The projections of $G_1$-translates of $B_1$ to the last $n-k$
coordinates $\R^{n-k}$ are bounded away from $\infty$.
\end{cl}

\begin{proof}
Since $U$ is contained in $B^\circ_2$, the last $n-k$ coordinates of
its complement $B_1$ are bounded away from $\infty$.
Moreover since
$\sum^{n}_{i=k+1}{|g(x)|_i^2}=\sum^{n}_{i=k+1}|x|_i^2$ for any $g
\in H$, by taking $t$ sufficiently large, we know that $g(B_1)\cap
U=\emptyset$. This means that the projections of $H$-translates of
$B_1$ to the last $n-k$ coordinates of $\R^{n-k}$ are bounded away from
$\infty$.

Now we consider general translates by elements of $G_1$.
Suppose, seeking a contradiction, that there is a sequence $\{g_k(B_1)\}$ of
distinct $G_1$-translates of $B_1$ whose projections to $\R^{n-k}$ go to $\infty$.
 Since $J$ stabilises $B_1$, we see that $g_k\in G_1-(H\cup J)$.

 On the other hand, since $U$ is a peak domain for $G_1$, it is precisely
 invariant under $H$ in $G_1$.
Take a point $y_0$ in $U$.
Since $g_k(y_0)$ is disjoint from $U$, the last $n-k$ coordinates of
$g_k(y_0)$ are bounded as $k \rightarrow \infty$.
Since $H$ acts on the first $k$-coordinates cocompactly, we can choose
$j_k \in H$ such that $j_kg_k(y_0)$ stays in a bounded set.

Since $j_k$ lies in $H$, we have $\sum_{i=k+1}^{n}(j_k(x))_i^2=\sum_{i=k+1}^{n}(x)_i^2$.
Therefore the projections of $j_kg_k(B_1)$ to $\R^{n-k}$ also go to $\infty$.
 Now Lemma \ref{key lemma}-(7) implies that $j_k g_k(y)\rightarrow \infty$ for all $y\in B_1$.
 By Lemma \ref{converge to limit points}, we see that, by choosing a subsequence
 if necessary, we may assume that $j_k g_k(y)\rightarrow \infty$ for all
 $y$ except for at most one point which is contained in the limit set of $G_1$.
 Since $y_0$ is contained in $U \subset \Omega(G_1)$, we have in
 particular that $j_kg_k(y_0) \rightarrow \infty$.
 This is a contradiction.
\end{proof}

Our claim shows that $U$ can be taken to be  disjoint from $T_1$.
Therefore, we have $U\subset C_1\cap B^\circ_2$. Since
$C_1\cap B^\circ_2$ is precisely invariant under $G_1$ in $G$, for
any $g\in G-G_1$, $g(U)\cap U=\emptyset$. Therefore, $U$ is a peak
domain at $x$ of $G$, which means that $x$ is a parabolic vertex for
$G$. Thus we have proved that all parabolic fixed points of $G$ are
parabolic vertices.

Next assume that $x$ is a limit point of $G$ which is not a
parabolic fixed point. Suppose that  $x$ is a translate of
a limit point $y$ of $G_m$. Since $y$ is not a parabolic fixed point
and $G_m$ is geometrically finite, by Proposition
\ref{limit set}, $y$ is a conical limit point of $G_m$,
hence also for $G$. If $x$ is not a translate of a limit point of
either $G_1$ or $G_2$, then by (10), it is a conical limit point for
$G$. Thus we have shown that any non-parabolic limit point of $G$ is
a conical limit point, and completed the proof of the ``if" part.

We shall now turn to show the ``only if"  part.
Assume that $G$ is geometrically finite.
Then $S$ is a strong $(J,G)$-block.
This implies that $B_m$ is a strong $(J,G_m)$-block for $m=1,2$ by (9).

Let $x$ be a parabolic fixed point of $G_1$. We assume that the rank
of $\mathrm{Stab}_{G_1}(x)$ is $k<n$, and shall prove that there is
a peak domain at $x$ for $G_1$. Since $B_1^\circ$ is contained in
$\Omega(G_1)$ by Lemma \ref{key lemma}-(3), $x$ cannot lie in
$G_1(B_1^\circ)$. Therefore,  $x$ lies in either  $G_1(S)$
or $B_2^\circ \cap C_1$. If $x\in G_1(S)$, then, since $B_1$ is a
strong $(J,G_1)$-block and $J$ is geometrically finite, there is a peak domain  at $x$ for $G_1$,
and we are done. If $x\in B^\circ_2\cap C_1$, then
$\mathrm{Stab}_G(x)=\mathrm{Stab}_{G_1}(x)$ since $B^\circ_2\cap
C_1$ is precisely invariant under $G_1$ in $G$ by Lemma \ref{C_1
invariant}. Therefore $\mathrm{Stab}_G(x)$ has rank $k < n$ in
particular. Since $G$ is geometrically finite, there is a peak
domain $U$  at $x$ for $G$, which is also a peak domain for $G_1$.

Now let $x$ be a limit point of $G_1$ which is not a parabolic fixed
point of $G_1$. We shall show that $x$ is a conical limit point of
$G_1$. Again we have only to consider the cases when $x \in G_1(S)$
and when $x \in B^\circ_2 \cap C_1$. If $x\in G_1(S)$, then there
are a point $y$ lying on $S$ and $g \in G_1$ such that
$x=g_1(y)$. Since $y$ lies on $\Lambda(J)$ by Lemma \ref{key
lemma}-(2), and $J$ is geometrically finite, it is a conical limit
point for $J$ by Proposition \ref{limit set}. This implies that $x$
is a conical limit point for $G_1$, and we are done in this case.

Suppose now that $x \in B^\circ_2 \cap C_1$.
 Since $B^\circ_2\cap C_1$ is precisely invariant under $G_1$, we have
 $\mathrm{Stab}_G(x)=\mathrm{Stab}_{G_1}(x)$.
Therefore $x$ is not a parabolic fixed point of $G$ either. Since
$G$ was assumed to be geometrically finite, $x$ is a conical limit
point for $G$ by Proposition \ref{limit set}. It follows from
Proposition \ref{conical limit} that  there is a sequence $\{h_k\}$
of distinct elements of $G$ such that $d(h_k(z),h_k(x))$ is bounded
away from zero for all $z\in \bar \R^n\backslash\{x\}$ and $h_k^{-1}(z_0)\rightarrow x$
for some $z_0\in \hyp^{n+1}$.

\begin{cl}
\label{distinct}
All the $h_k(S)$ are distinct passing to a subsequence if necessary.
\end{cl}

\begin{proof}
Recall that we assumed that  $S$ is precisely invariant under $J$ in both $G_1$ and $G_2$.
Therefore, $S$ is precisely invariant under $J$ in $G$ by Lemma
\ref{precisely embedded}.

Now, suppose, seeking a contradiction,  that all the $h_k(S)$ are
the same after passing to a subsequence. Then $h^{-1}_1\circ
h_k(S)=S$ for every $k$. Therefore for each $k$, there is an element
$j_k\in J$ such that $h_k=h_1\circ j_k$, with $j_1=id$. Since $h_k$
are distinct elements of $G$, all $j_k$ are distinct and by Lemma
\ref{converge to limit points} and the fact that $d(h_k(z),h_k(x))$
is bounded away from $0$ for all $x \neq z$, we may assume that
there are two distinct points $x', z'$ such that $h_k(z)=h_1\circ
j_k(z)\rightarrow z'$ for all $z\in \bar \R^n\setminus
\{x\}$ and $h_k(x)=h_1\circ j_k(x)\rightarrow x'$. This
implies that $j_k(z)\rightarrow h^{-1}_1(z')$ for all $z\in \bar
\R^n\setminus \{x\}$ and $j_k(x)\rightarrow
h^{-1}_1(x')$. Since $z'\neq x'$ and
$j_k^{-1}(z_0')\rightarrow x$, where $z_0'=h_1^{-1}(z_0)\in
\hyp^{n+1}$, it follows that $x$ is a conical limit point of $J$ by
Proposition \ref{conical limit}. Since we assumed that $x\in C_1\cap
B_2^\circ$, we have $x\notin \Lambda(J)$. This
is a contradiction and we have completed the proof of Claim
\ref{distinct}.
\end{proof}


\begin{cl}
\label{positive}
By taking a subsequence we can assume $h_k>0$ for all $k$.
\end{cl}

\begin{proof}
Suppose, on the contrary, that  $h_k<0$ for all $k$ after passing to a subsequence.
Since all $h_k(S)$ are distinct, the $h_k$ belong to distinct cosets of $J$ in $G$.
By (5), we have $\mathrm{diam}(h_k(S)) \rightarrow 0$.
Since we assumed $h_k <0$, the set $h_k(B_2)$ cannot contain $S$ inside, hence
is contained in the smaller part of $\bar \R^{n} \setminus h_k(S)$.
Therefore, we have $\mathrm{diam}(h_k(B_2)) \rightarrow 0$.
Recall that we are considering the case when $x \in B^\circ_2 \cap C_1$.
Therefore, we have $d(h_k(z), h_k(x)) \rightarrow 0$ for all $z \in B_2$.
This contradicts the fact that $d(h_k(z), h_k(x))$ is bounded away from $0$ for $z \in \bar \R^n \setminus \{x\}$.
 Thus we have completed the proof of Claim \ref{positive}.
\end{proof}

Now we return to the proof of (11). Note that we have only to
consider the case when $h_k$ is not contained in $G_1$, for
otherwise $x$ is a conical limit point of $G_1$ by Proposition
\ref{conical limit}. Therefore, we can assume that
$|h_k|>1$. Express $h_k$ in a normal form
$h_k=\gamma_{k_l}\circ\cdots\circ \gamma_{k_1}$. Set $g_k = h_k
\circ \gamma_{k_1}^{-1}$. Then $g_k$ is negative.

First consider the case when  $g_k=g\circ j_k$ for some $g\in G$
with some $j_k\in J$. Then $d(h_k(z),h_k(x))=d(g\circ j_k\circ
\gamma_{k_1}(z),g\circ j_k\circ \gamma_{k_1}(x))$. By Lemma
\ref{converge to limit points}, we may assume that there are two
distinct points $x',z'$ such that $g\circ j_k\circ
\gamma_{k_1}(z)\rightarrow z'$ for all $z\in \bar \R^n \setminus
\{x\}$ and $g\circ j_k\circ \gamma_{k_1}(x)\rightarrow x'$. It
follows that $j_k\circ \gamma_{k_1}(z)\rightarrow g^{-1}(z')$ for
all $z\in \bar \R^n\setminus \{ x\}$, $j_k\circ
\gamma_{k_1}(x)\rightarrow g^{-1}(x')$ and
$(j_k\circ\gamma_k)^{-1}(g^{-1}(z_0))\rightarrow x$, where
$g^{-1}(z_0)\in \hyp^{n+1}$. It follows from Proposition
\ref{conical limit} that $x$ is a conical limit point of $G_1$.

 Suppose next that $g_k$ is not expressed as $g \circ j_k$.
 Then by Claim \ref{distinct},  $g_k(S)$ are all distinct.
 Applying the proof of Claim \ref{positive} to $g_k$, we see that
 $\mathrm{diam}(g_k(B_2)) \rightarrow 0$.
 Now, $Q_1=G_1(B_1^\circ)$ is invariant under $G_1$, hence so is its complement $R_1$.
It follows that $h_k(R_1)=g_k(R_1)$. Since $R_1$ is contained in
$B_2$, we have
$\mathrm{diam}(h_k(R_1))=\mathrm{diam}(g_k(R_1))\rightarrow 0$. By
the same argument as in the proof of Claim \ref{positive} and the
fact that $S\subset R_1$, this is a
contradiction. Thus we have completed the proof of (11).
\end{proof}

\begin{rem}
The condition that $(B^\circ_1, B^\circ_2)$ is a proper interactive
pair in Theorem \ref{main} is necessary, as the following example shows.
\end{rem}

\begin{exa}
Set
$$J=\langle \left(\begin{array}{cc}
                      1 & 1 \\
                      0 & 1
         \end{array}  \right),\,\  \left(
         \begin{array}{cc}
                      0 & 1 \\
                      -1 & 0
         \end{array}
    \right)\rangle,\;\; g_1= \left(\begin{array}{cc}
                      i & 0 \\
                      0 & -i
         \end{array}  \right),\;\; g_2=\left(\begin{array}{cc}
                      0 & i \\
                      i & 0
         \end{array}  \right)
$$ and

$$G_1=\langle J,g_1\rangle,
\;\; G_2=\langle J,g_2\rangle.$$

We use the following symbols:

$$S=\{x\in \bar \R^2:\; x_2=0\},\;\; B_1=\{x\in \bar \R^2:\; x_2\leq 0\}
\;\;\mbox{and}\;\;B_2=\{x\in \bar \R^2:\; x_2\geq 0\}.$$

Then the following hold.

\begin{enumerate}
\item $J$ is geometrically finite.

\item $S=\Lambda(J)=\Lambda(G_1)=\Lambda(G_2)$.

\item $G_1=J\cup g_1J$ and $G_2=J\cup g_2J$.

\item  Each $B_m$ is a $(J,G_m)$-block for $m=1,2$.

\item $(B^\circ_1,B^\circ_2)$ is an interactive pair, but
$(B^\circ_1,B^\circ_2)$ is not proper.

\item  $G\not=G_1*_JG_2$.
\end{enumerate}
\end{exa}

The assertion (1) is obvious since $J$ is a finitely generated
Fuchsian group.
To prove (2), set $\displaystyle w=\frac{p}{r}$, where $p$ and $r$ are
integers and $r\neq0$, and $j=\left(
         \begin{array}{cc}
                      1-pr & p^2 \\
                      -r^2 & 1+pr
         \end{array}
    \right)$.
Then $j\in J$ is a parabolic element having $w$ as its fixed point.
Therefore, every rational number is a parabolic fixed point of $J$.
Now (2) follows from Lemma 5.3.3 in \cite{Bea}.
The proofs of $(3)$, $(4)$ and $(5)$ are trivial.
We can verify
(6) by checking that for  a $(1,2)$-form $g_1g_2g_1g_2$, we have
 $\Phi(g_1g_2g_1g_2)=id$.

\section{Applications}
Following  \cite{Wan} or \cite{Wat}, we denote by $PSL(2,\Gamma_n)$ the
$n$-dimensional Clifford matrix group. Then $PSL(2,\Gamma_n)$ is
isomorphic to $M(\bar \R^n)$ (cf. \cite{Ahl}).

\begin{exa}
\label{first example}
Let $n\geq 4$, and set

 $$j= \left(\begin{array}{cc}
                      e_1e_2 & 0 \\
                      0 & e_1e_2
         \end{array}  \right),
         \ g_1= \left(
         \begin{array}{cc}
                      0 & e_{n-1} \\
                      e_{n-1} & 0
         \end{array}
    \right)\;\;\mbox{and}\;\;
    \ g_2= \left(
         \begin{array}{cc}
                      0 & 2e_{n-1} \\
                      \frac{1}{2}e_{n-1} & 0
         \end{array}
    \right).$$
We also set
$$J=<j>,\;\;G_1=<j, g_1>\;\;\mbox{ and}\;\; G_2=<j,g_2>.$$
\end{exa}

Since $J$ is a finite group, it is geometrically finite. Set
$S=\{x\in \bar \R^n:\; |x|=\sqrt{2}\},\;\; B_1=\{x\in \bar \R^n:\;
|x|\geq\sqrt{2}\}\;\; \mbox{and}\;\; B_2=\{x\in \bar \R^n:\;
|x|\leq\sqrt{2}\}.$ Obviously,  $G_m=\{id, j, g_m, jg_m\}$ is
geometrically finite for $m=1,2$. Set
$x=x_0+x_1e_1+\cdots+x_{n-1}e_{n-1}+x_ne_n\in \hyp^{n+1}$. Then we
have $\mathrm{Fix}(\widetilde{j})=\{x\in \hyp^{n+1}\cup \R^n:\;
x_1=x_2=0\}\cup\{\infty\},$ $\mathrm{Fix}(\widetilde{g_1})=\{x\in
\hyp^{n+1}\cup \R^n:\; |x|=1\;\
       \mbox{and}\;\
       x_{n-1}=0\},$
$\mathrm{Fix}(\widetilde{g_2})=\{x\in \hyp^{n+1}\cup \R^n:\;
|x|=2\;\
       \mbox{and}\;\ x_{n-1}=0\}$
and $\mathrm{Fix}(\widetilde{jg_m})=\{x\in \hyp^{n+1}\cup \R^n:\;
|x|=m\;\
       \mbox{and}\;\ x_1=x_2=x_{n-1}=0\}$, where $\widetilde{h}$
       denotes the Poincar{\'e} extension of $h\in M(\bar \R^n)$ in
       $\B^{n+1}$
and $\mathrm{Fix}(\widetilde{h})=\{x\in \B^{n+1}:\;
\widetilde{h}(x)=x\}$. Therefore for each $m$ ($m=1,2$),
$\mathrm{Fix}(\widetilde{jg_m})= \mathrm{Fix}(\widetilde{j})\cap
\mathrm{Fix}(\widetilde{g_m})$.

We put $a=e_1+e_n$. It is obvious that $a$ is not fixed by any
nontrivial element in either $G_1$ or $G_2$. For any non-trivial
element $h\in M(\bar \R^n)$, if we set $H_h=\{x\in \hyp^{n+1}:\;
d_h(x,a)\leq d_h(x, ha)\}$, then $H_j=\{x\in \hyp^{n+1}:\; x_1\geq
0\},\;\ H_{g_1}=\{x\in \hyp^{n+1}:\; |x|\geq 1\},$ $H_{jg_1}=\{x\in
\hyp^{n+1}:\; |x+2e_1|\geq\sqrt{5}\},\;\ H_{g_2}=\{x\in
\hyp^{n+1}:\; |x|\leq 2\},$ and $H_{jg_2}=\{x\in \hyp^{n+1}:\;
|x-4e_1|\leq 2\sqrt{5}\}.$

For each $m$ ($m=1,2$), set $P_m=H_j\cap H_{g_m}\cap H_{jg_m}$.
 Then $P_m$ is the closure of the Dirichlet fundamental polyhedron
centered at $a$ for $G_m$ in $\hyp^{n+1}$ (c.f. \cite{Rat}) and
$P_m=H_j\cap H_{g_m}$. We consider $ D_1=\{x\in \R^n:\; x_1\geq
0,\;\ |x|\geq 1\}\setminus (\{x\in \R^n:\; |x|=1\;\
       \mbox{and}\;\
       x_{n-1}\leq0\}\cup\{x\in \R^{n}:\; x_1=0, x_2\leq0\}),$
$D_2=\{x\in \R^n:\; x_1\geq0,\;\ |x|\leq2\}\setminus(\{x\in \R^n:\;
|x|=2\;\
       \mbox{and}\;\
       x_{n-1}\leq0\}\cup\{x\in \R^{n}:\; x_1=0, x_2\leq0\}).
       $

It is easy to see that for each $m$ ($m=1,2$), $B_m$ is a $(J,
G_m)$-block and precisely invariant under $J$ in $G_m$, that $D_m$ is a
fundamental set for $G_m$ with $J(D_m\cap B_m)=B_m\cap$
$^\circ\Omega(J)$, that $(D_m\cap B_{3-m})^\circ\not=\emptyset$, and
that $D_1\cap S=D_2\cap S$.
Since $(B^\circ_1, B^\circ_2)$ is a proper interactive pair, $<G_1,
G_2>=G_1*_JG_2$.
It is obvious that $G_1$ and $G_2$ are
geometrically finite. It follows from Theorem \ref{main} that $G$
is also geometrically finite.

In this example, the amalgamated free product $G_1*_JG_2$ is
elementary.
The following two examples give non-elementary groups.

\begin{exa}
Suppose that  $n\geq 3$ and let $j$, $J$ be the same as in Example \ref{first example}.
We set
$$g_1= \left(
         \begin{array}{cc}
                      1 & 0 \\
                      2 & 1
         \end{array}
    \right)\;\;\mbox{and}\;\;
    \ g_2= \left(
         \begin{array}{cc}
                      1 & 5 \\
                      0 & 1
         \end{array}
    \right),$$
$G_1=<j, g_1>$, $G_2=<j,g_2>$, and $S=\{x\in \bar \R^n:\; |x|=2\}$,
$B_1=\{x\in \bar \R^n:\; |x|\geq 2\}$, $B_2=\{x\in \bar \R^n:\;
|x|\leq2\}$.
\end{exa}

We define two domains by $\displaystyle D_1=\{x\in \R^n:\;
x_1\geq0,\;\ |x+\frac{1}{2}|\geq \frac{1}{2},\;\
|x-\frac{1}{2}|>\frac{1}{2}\}\setminus \{x\in \R^n:\; x_1=0,
x_2\leq0\}$ and $\displaystyle D_2=\{x\in \R^n:\; x_1\geq0,\;\
\frac{-5}{2}< x_0\leq\frac{5}{2}\}\setminus \{x\in \R^n:\; x_1=0,
x_2\leq0\}.$ Then  the discussion similar to the one in Example
\ref{first example} shows that for each $m$ ($m=1,2$),
\begin{enumerate}
\item $B_m$ is a ($J, G_m$)-block and precisely invariant under
$J$ in $G_m$;
\item $D_m$ is a fundamental set for $G_m$ satisfying that
$J(D_m\cap B_m)=B_m\cap $$^\circ\Omega(J)$;
\item $(B_1^\circ, B_2^\circ)$ is a proper interactive pair;
\item  $(D_m\cap B_{3-m})^\circ\not=\emptyset$; and
\item  $D_1\cap S=D_2\cap S$.
\end{enumerate}

Theorem \ref{main} implies that
\begin{enumerate}
\item $G=<G_1, G_2>=G_1*_JG_2$;

\item $G$ is geometrically finite since both $G_1$ and $G_2$
are geometrically finite.
\end{enumerate}

\begin{exa}
Suppose that  $n\geq 5$, and let $j$, $J$ be the same as in Example \ref{first example}.
We set
$$g_1= \left(
         \begin{array}{cc}
                      1 & 0 \\
                      2e_3 & 1
         \end{array}
    \right)\;\;\mbox{and}\;\;
    \ g_2= \left(
         \begin{array}{cc}
                      1 & 5e_{n-1} \\
                      0 & 1
         \end{array}
    \right),$$
and $G_1=<j, g_1>$, $G_2=<j,g_2>$. We define $S, B_1, B_2$ by
$S=\{x\in \bar \R^n:\; |x|=2\}$, $B_1=\{x\in \bar \R^n:\;|x|\geq
2\}$, and $B_2=\{x\in \bar \R^n:\; |x|\leq2\}$.
\end{exa}

We define two domains $D_1, D_2$  by $\displaystyle D_1=\{x\in
\R^n:\; x_1\geq0,\;\ |x+\frac{e_3}{2}|\geq \frac{1}{2} ,\;\
|x-\frac{e_3}{2}|>\frac{1}{2}\}\setminus\{x\in \R^n:\; x_1=0,
x_2\leq0\}$ and $D_2=\{x\in \R^n:\; x_1\geq0,\;\ -\frac{5}{2}<
x_{n-1}\leq\frac{5}{2}\}\setminus \{x\in \R^n:\; x_1=0, x_2\leq0\}.$
Then we again have the following.

\begin{enumerate}
\item $B_m$ is a ($J, G_m$)-block and precisely invariant under
$J$ in $G_m$;

\item $D_m$ is a fundamental set for $G_m$ satisfying that
$J(D_m\cap B_m)=B_m\cap $$^\circ\Omega(J)$;

\item $(B_1^\circ, B_2^\circ)$ is a proper interactive pair;

\item $(D_m\cap B_{3-m})^\circ\not=\emptyset$;

\item $D_1\cap S=D_2\cap S$.
\end{enumerate}

Therfore, Theorem \ref{main} implies that

\begin{enumerate}
\item  $G=<G_1, G_2>=G_1*_JG_2$; and
\item $G$ is geometrically finite since both $G_1$ and $G_2$
are geometrically finite.
\end{enumerate}

\makeatletter
\newlength{\arclength}
\newcommand{\arc}[1]{\settowidth{\arclength}{#1}%
                     \setlength{\unitlength}{0.01\arclength}%
                     \providecommand{\@arc}{}%
                     \renewcommand{\@arc}{%
                     \begin{picture}(100,30)%
                     \qbezier(0,0)(50,35)(100,0)%
                     \end{picture}}%
                     \stackrel{\@arc}{#1}}
                     \makeatother

$\arc{ab}$
\end{document}